\documentclass[11pt]{article}

\usepackage[margin=.7in]{geometry}
\usepackage{amsfonts}
\usepackage[mathcal]{euscript}
\usepackage{amsmath,amssymb}
\usepackage{amsthm}
\usepackage{graphicx}
\usepackage{pgfplots}
\pgfplotsset{compat=1.14}

\usepackage{caption}
\usepackage{subcaption}

\usepackage{hyperref}

\usepackage{xcolor}

\usepackage[ruled,linesnumbered]{algorithm2e}

\usepackage{tabularx}
\usepackage[normalem]{ulem}

\usepackage[affil-it]{authblk}
\usepackage{stmaryrd}
\usepackage{mathtools}
\usepackage{bm}
\usepackage{array}
\usepackage{empheq}

\newcommand{\mi}{\textrm i}

\newcommand{\pol}[1]{\mathbb{P}_{#1}} 
\newcommand{\hompol}[1]{\tilde{\mathbb{P}}_{#1}} 

\newcommand{\homharm}[1]{\tilde{\mathbb H}_{#1}}

\newtheorem{thm}{Theorem}
\newtheorem{cor}{Corollary}
\newtheorem{prop}{Proposition}
\newtheorem{dfn}{Definition}
\newtheorem{rmk}{Remark}
\newtheorem{lmm}{Lemma}

\newcommand{\modTC}[1]{\textcolor{magenta}{#1}}

\let\oldnl\nl
\newcommand{\nonl}{\renewcommand{\nl}{\let\nl\oldnl}}

\DeclareCaptionLabelFormat{andtable}{#1~#2  \&  \tablename~\thetable}

\usepackage{tikz}

\tikzstyle{arrow} = [thick,->,>=stealth]

\pgfarrowsdeclare{mylatex'}{mylatex'}
{
	\newdimen\len
	\len=\pgfgetarrowoptions{mylatex'}
	\pgfarrowsleftextend{-0.4\len}
	\pgfarrowsrightextend{0.6\len}
}
{
	\newdimen\len
	\len=\pgfgetarrowoptions{mylatex'}
	\pgfpathmoveto{\pgfqpoint{0.6\len}{0\len}}
	\pgfpathcurveto
	{\pgfqpoint{0.35\len}{0.05\len}}
	{\pgfqpoint{-0.1\len}{0.15\len}}
	{\pgfqpoint{-0.4\len}{0.375\len}}
	\pgfpathcurveto
	{\pgfqpoint{-0.15\len}{0.1\len}}
	{\pgfqpoint{-0.15\len}{-0.1\len}}
	{\pgfqpoint{-0.4\len}{-0.375\len}}
	\pgfpathcurveto
	{\pgfqpoint{-0.1\len}{-0.15\len}}
	{\pgfqpoint{0.35\len}{-0.05\len}}
	{\pgfqpoint{0.6\len}{0\len}}
	\pgfusepathqfill
}
\pgfsetarrowoptions{mylatex'}{8pt}
\pgfkeys{/tiplen/.default=8pt, /tiplen/.code={\pgfsetarrowoptions{mylatex'}{#1}}}

\title{
Quasi-Trefftz spaces for a second order time-harmonic Maxwell's equation
}

\author{Ilaria Fontana and Lise-Marie Imbert-G\'erard (University of Arizona)}

\begin{document}
\maketitle
\tableofcontents

{\bf Acknowledgments}\\
This material is based upon work supported by the U.S. Department of Energy, Office
of Science, Office of Advanced Scientific Computing Research, under Award Number DE-SC0024246.
\newpage

{\bf Abstract}\\
Quasi-Trefftz methods are a family of Discontinuous Galerkin methods relying on equation-dependent function spaces.
This work is the first study of the notion of local Taylor-based polynomial quasi-Trefftz spaces for a system of Partial Differential Equations (PDEs). 
These discrete spaces are meant to discretize a quasi-Trefftz weak formulation, and are introduced here for electro-magnetic wave propagation in inhomogeneous media, governed by a second order formulation of Maxwell's equation with variable coefficients.
Thanks to an adequate Helmholtz decomposition for spaces of homogeneous polynomial vector fields, 
the outcome is the proposed quasi-Trefftz spaces with same best approximation properties as the corresponding standard polynomial spaces but with an optimal dimension, as well as a procedure to construct a quasi-Trefftz basis. 

\section{Introduction}
In the field of numerical partial differential equations (PDEs), there exist various methods leveraging problem-dependent basis functions \cite{BOTTASSO20023391,DEMKOWICZ20101558,Altmann_Henning_Peterseim_2021,Peterseim2014EliminatingTP,MELENK1996289}.
Specifically for problems modeling wave propagation phenomena, different high-order methods incorporate the expected oscillating behavior of the solution to the PDE problem within the discrete basis \cite{FARHAT20016455,DECKERS20111117,DesmetThesis}.

Among them, a specific type of discontinuous Galerkin methods, referred to as Trefftz methods, fundamentally rely on local spaces of exact solutions to the governing PDE 
 \cite{BrunoCRAS,doi:10.1137/S0036142995285873,doi:10.1137/S1064827503422233,GABARD20071961,https://doi.org/10.1002/nme.4961,PWHipt++,HMPh}. For convenience, in this work, the notion of Trefftz property for a function refers to being a local exact solution to the governing PDE. A weak formulation holding only for Trefftz functions is then discretized on local discrete Trefftz spaces.
These methods have been developed in particular for time-harmonic wave propagation problems, mainly thanks to discrete spaces of local Plane Waves (PWs), but also spherical waves or Bessel-based waves (generalized harmonic polynomials) \cite{Luostari_Huttunen_Monk_2012}.
From the point of view of the size of the discrete system, such methods may need much smaller discrete spaces in order to reach the same order of accuracy compared to standard Galerkin methods, resulting in a considerably smaller linear system.
From the point of view of the assembly of the matrix, such methods present two main advantages. 
First, thanks to the Trefftz property, the weak formulation contains only integrals on the boundary of mesh elements, hence the domains of integration are lower-dimensional compared to standard Galerkin method.
Second, there are closed formulas for integrating products of PWs and their derivatives.
Both facts reduce considerably to cost of assembly of the matrix.
Surveys of some of these methods can be found in \cite{Hiptmair2016,DECKERS2014550}, and some related computational aspects are covered in \cite{LIEU2016105,HUTTUNEN200227}.

Trefftz-like methods also have known limitations. 
From the point of view of wave-like bases, there exists a conditioning problem that can be simply illustrated considering the difference between two PWs propagating respectively in the directions $ \mathbf d_1$ and $ \mathbf d_2$: 
\begin{equation*}
e^{\mi\kappa \mathbf d_1 \cdot\mathbf x} -e^{\mi\kappa \mathbf d_2 \cdot\mathbf x}= 2\mi \exp\left( \frac{\mi\kappa(\mathbf d_1+\mathbf d_2)\cdot\mathbf x}{2} \right)\sin \frac{\kappa(\mathbf d_1-\mathbf d_2)\cdot\mathbf x}{2}
\rightarrow 0 
\textbf{ as } |\mathbf d_1-\mathbf d_2|\rightarrow 0,
\end{equation*}
as the linear independence of these two PWs deteriorates as their directions of propagation get closer.
Moreover, the range of applications of Trefftz-like methods is limited by the need for discrete spaces of exact solutions to the governing PDEs with good approximation properties and explicit basis functions.
For most problems beyond wave propagation in homogeneous media, 
 such discrete spaces are not available.

A class of so-called quasi-Trefftz methods was introduced to circumvent this limitation.
These methods relax the Trefftz property into a so-called quasi-Trefftz property: being a local approximate solution to the governing PDE, by imposing that a Taylor polynomial of the differential operator applied to the function vanishes up to a desired order.  
These methods were first developed for the Helmholtz equation {with} variable wave number based on so-called Generalized { PW} basis functions \cite{10.1093/imanum/drt030,LMIGinterp,LMP}, and later for the wave equation \cite{IGMS}, the Shr\"odinger equation \cite{GM}, and more generally to classes of elliptic problems \cite{10.1093/imanum/drae094} and elliptic-hyperbolic problems \cite{10.1051/m2an/2026081}. The methods have so far not been studied for problems governed by systems of PDEs.
The development of quasi-Trefftz methods for problems governed by systems of PDEs is of great interest to enable accurate and efficient simulations in electromagnetics and many other applications. In this context, this work is the first study of Taylor-based polynomial quasi-Trefftz discrete spaces for a problem governed by a system of PDEs. As in previous works, such spaces are {\bf subspaces standard polynomial spaces that enjoy the same high-order best approximation properties but have a significantly smaller dimension}.
For problems governed by a system of PDEs, the coupling between the different components of the unknown naturally raises new questions compared to the scalar case. In this work, the divergence condition plays a central role in defining the smallest possible quasi-Trefftz spaces.

Thanks to the tools developed in the present work, Taylor-based polynomial quasi-Trefftz discrete spaces were later studied for other PDE systems, including a first-order formulation of Helmholtz equation \cite{AndL1Helm} and for a first-order formulation of Maxwell's equations \cite{fontana2026quasitrefftzspacesfirstordertimeharmonic}.

\subsection{Maxwell's equation}
This article focuses on Taylor-based polynomial quasi-Trefftz spaces for a second-order time-harmonic form of Maxwell's equations.
Given a smoothly varying coefficient $\varepsilon$, consider
\begin{equation}\label{eq:Max}
\nabla\times\nabla\times\mathbf V -\varepsilon \mathbf V = \mathbf 0.
\end{equation}
Since $\nabla\cdot (\nabla \times\mathbf V) = \mathbf 0$, any solution $\mathbf V$ to Equation \eqref{eq:Max} necessarily satisfies the so-called divergence condition $\nabla\cdot(\varepsilon \mathbf V)=0$.


\begin{rmk}
In the present work, the variable coefficient $\varepsilon$ is assumed to be scalar, yet a matrix-valued dielectric tensor $\varepsilon$ is key to the studying the interaction of various propagation modes, for example in the cold plasma model for radio-frequency wave propagation in a magnetized fusions plasma.
Extending this work to the case of a matrix-valued variable coefficient $\varepsilon$ is straightforward, and simply requires  $\varepsilon$ to be smooth component-wise and for its eigenvalues to be bounded away from zero, with no restriction on their sign. 
Moreover, for a matrix-valued coefficient, the product rule for the divergence 
reads 
 $$
 \nabla\cdot(\varepsilon\mathbf V)
 = \sum_{g'}\sum_{g}\partial_g(\varepsilon_{gg'})V_{g'}+\sum_{g}\sum_{g'}\varepsilon_{gg'}\partial_g(V_{g'}).
 $$
 \end{rmk}

There exists a rich literature dedicated to the mathematical theory of Maxwell's equations and their applications, see for instance \cite{Mullbook69,doi:10.1142/2938,nedelec2001acoustic,rodriguez2010eddy,garrity2015electricity,MAFEboook,colton2019inverse,donnevert2020maxwell,maxworth2024one}. 
The divergence condition and the Helmholtz decomposition of vector fields as the sum of a divergence-free (or solenoidal) component and a curl-free (or irrotational) component play a fundamental role in this context \cite{vivette,6365629}, despite the existence of many different formulations of these equations.

Many numerical methods and tools have also been developed specifically for systems of PDEs, and in particular for Maxwell's equations \cite{Inan_Marshall_2011,doi:10.1137/110830320,doi:10.1137/1.9781611975543,10.1007/s00211-017-0939-x,Cotter_2023}.
These for instance include the Finite Difference Yee scheme, the Finite Element Exterior Calculus, and various types of wave propagation algorithms for hyperbolic systems.
Here again, the divergence condition remains fundamental, and
the choice of discrete spaces is crucial for Galerkin methods, for Maxwell's equations  \cite{10.1093/acprof:oso/9780198508885.001.0001,doi:10.1137/S003614290241790X,10.1002/lpor.201000045,Xie_Wang_Zhang_2013,SHI2018147,doi:10.1137/1.9781611975543}. 
For example, N\'ed\'elec finite elements can be used to discretize $H(curl)$ and Raviart-Thomas finite elements can be used to discretize $H(div)$, and more generally Whitney forms can be used to respect the relationship between the gradient, divergence and curl operators \cite{ArnoldFalkWinther2006,Nedelec1980,Hiptmair_2002,Bossavit1998,Bossavit1988Whitney,Whitney1957}. 
The dimension of such spaces is also of interest \cite{doi:10.1137/22M1544579}.
More generally, exact sequences of finite element spaces have been developed to discretize differential complexes, such as the deRham complex for Maxwell problems \cite{DEMKOWICZ200029,Hiptmair_2002,10.1007/0-387-38034-5_2,10.1007/s00211-017-0939-x}, in particular for stability \cite{Arnold02DCnNS}. 
Many versions of discrete Helmholtz decompositions have also been proposed, within the Galerkin framework and beyond \cite{cite-key,discrVFdecomp,discreteHD2025}.

In the case of  Trefftz-Discontinuous Galerkin (DG) methods for wave propagation in homogeneous media, the discrete space of plane waves (PWs) includes a compatibility condition between the propagation directions and the polarizations to guarantee that the PWs satisfy the governing PDE \cite{AndMax}.
In this first work on quasi-Trefftz spaces for Maxwell's equation, an adequate {\bf Helmholtz decomposition of homogeneous polynomial vector fields}, including a harmonic component, is proposed. This is the {\bf key to the study of polynomial quasi-Trefftz spaces for Maxwell's equation \eqref{eq:Max}}.

\subsection{Notation}
Throughout the article, 
while $\mathbb N$ refers to the set of positive integers, $\mathbb N_0$ refers to the set of non-negative integers;
the kernel of an operator is denoted $\ker$ and range $\mathcal R$.
Moreover, 
$g$ refers to an integer,  
 the dimension of the ambient space is equal to $3$. Then $\mathbf x=(x_g,g\text{ from } 1 \text{ to } 3)\in\mathbb R^3$ refers to points in the ambient space, 
$\mathbf i=(i_g,g\text{ from } 1 \text{ to } 3)\in(\mathbb N_0)^3$ refers to a multi-index {with length $\lvert\mathbf i\rvert \coloneqq i_1+i_2+i_3$}, 
partial derivatives are denoted $\partial_{x_g}$ and $\partial^{\mathbf i}=\partial_{x_1}^{i_1}\partial_{x_2}^{i_2}\partial_{x_3}^{i_3} $, 
and finally the canonical basis of $\mathbb R^3$ is denoted $\{\mathbf e_g, 1\leq g\leq 3\}$.
The spaces of scalar and vector-valued polynomials of degree at most equal to $p\in\mathbb N_0$ are denoted respectively $\mathbb P_p$ and $(\mathbb P_p)^3$. 
The scalar polynomial zero is denoted $\mathbf 0$, the vector-valued polynomial zero is denoted $\mathbf 0^3$, and monomials are denoted $\displaystyle \forall \mathbf i \in(\mathbb N_0)^3, \mathbf X^{\mathbf i}=\prod_{g=1}^3 (X_{g})^{i_g}$.


The work presented here is limited to {\bf local} properties. 
More precisely, we only consider quasi-Trefftz spaces defined in the neighborhood of a point $\mathbf x_0\in\mathbb R^d$, so $\mathbf x_0$ will be omitted in the following notation. Spaces of scalar homogeneous polynomials are denoted 
$$
\forall l\in\mathbb N_0,\
\widetilde{\mathbb P}_l = Span\left\{ (\mathbf X-\mathbf x_0)^{\mathbf i}, \mathbf i \in(\mathbb N_0)^3, |\mathbf i|=l\right\}
$$
The graded structure of polynomial spaces can then be described as
$$
\forall l\in\mathbb N_0,\
 \mathbb P_p = \bigoplus_{l=0}^p \widetilde{\mathbb P}_l
\text{ and }
 (\mathbb P_p)^3 = \bigoplus_{l=0}^p \left(\widetilde{\mathbb P}_l\right)^3.
 $$
\begin{rmk}
As a reminder, dimensions of polynomial spaces in three variables are given by
$$
\forall p\in\mathbb N_0,\
\dim\widetilde{\mathbb P}_p = \frac 12(p+1)(p+2)
\text{ and }
\dim\mathbb P_p = \frac 16 (p+1)(p+2)(p+3).
$$
\end{rmk}
For $n\in\mathbb N_0$, the space of functions with continuous derivatives up to order $n$ in a neighborhood of $\mathbf x_0$ is denoted $\mathcal C^n(\mathbf x_0)$.
{
Given a vector field $\mathbf W\in (\mathcal C^{p+1}(\mathbf{x}_0))^3$, we denote by $\lVert \mathbf W \rVert_{\infty, \mathcal B(\mathbf{x}_0, r)}$ its sup-norm and by $\lvert\mathbf W\rvert_{\mathcal C^{p+1}(\mathcal B(\mathbf{x}_0, r))}$ its $\mathcal C^{p+1}$ seminorm over the ball centered in $\mathbf x_0$ and radius $r$, that is
\begin{equation*}
    \lVert \mathbf W \rVert_{\infty, \mathcal B(\mathbf{x}_0, r)} \coloneqq \sup_{\lVert\mathbf x - \mathbf x_0\rVert \leq r} \lVert \mathbf W(\mathbf x) \rVert
    \qquad\text{and}\qquad
    \lvert\mathbf W\rvert_{\mathcal C^{p+1}(\mathcal B(\mathbf{x}_0, r))} \coloneqq \max_{\mathbf i\in (\mathbb N_0)^3,\, \lvert\mathbf i\rvert = p+1} \lVert \mathbf \partial^{\mathbf i} \mathbf W \rVert_{\infty, \mathcal B(\mathbf{x}_0, r)}.
\end{equation*}
}
For a $\mathcal C^k(\mathbf x_0)$ function $\varphi$, its $k$th Taylor polynomial at $\mathbf x_0$, denoted $T_k[\varphi] \in\mathbb P_k$, is given by
$$
T_k[\varphi] =\sum_{l=0}^k \sum_{|\mathbf i| = l } \frac 1{\mathbf i!} \partial^{\mathbf i} \varphi(\mathbf x_0) (\mathbf X-\mathbf x_0)^{\mathbf i}.
$$
For convenience, for a $\varphi\in\mathcal C^k(\mathbf x_0)$, then $\varphi_k\in\widetilde{\mathbb P}_k$ will be used to denote the homogeneous component of degree $k$ of $T_k[\varphi] \in\mathbb P_k$.

The general concept of Taylor-based quasi-Trefftz spaces relies on the following notion of quasi-Trefftz property.
\begin{dfn}\label{def:atprop}
Assume $q\in\mathbb N_0$ and $d\in\mathbb N$.
Given a linear partial differential operator $\mathcal L$ with smooth variable coefficients in $\mathcal C^{q}(\mathbf x_0)$, the Taylor-based quasi-Trefftz property of order $q$ for a function $\varphi$ is
$$
T_q[\mathcal L  \varphi] = \mathbf 0^d
$$
{where $d$ is the dimension of the image of $\mathcal L$.}
\end{dfn}

\subsection{Quasi-Trefftz spaces: scalar equations versus systems of equations}
Quasi-Trefftz spaces of GPWs for various scalar equations were the focus of \cite{LMIGinterp,roadmap,doi:10.1137/20M136791X}, including a polynomial space in the latter.
In \cite{imbertgerard2025localtaylorbasedpolynomialquasitrefftz}, a general framework was proposed to study polynomial quasi-Trefftz spaces in any dimension for differential operators of order $\gamma\in\mathbb N$ of the form $\displaystyle \mathcal L := \sum_{m=0}^\gamma \sum_{ |\mathbf j|=m} c_{\mathbf j}\left( {\mathbf x}  \right) \partial_{{\mathbf x}}^{\mathbf j}$  with smooth variable coefficients in the vicinity of a point $\mathbf x_0$ with at least one $\mathbf j\in(\mathbb N_0)^3$ such that $|\mathbf j|=\gamma$ and $c_{\mathbf j}(\mathbf x_0)\neq 0$. There, the quasi-Trefftz spaces and the quasi-Trefftz operators are defined as
$$
\mathbb Q\mathbb T_p = \{Q\in \mathbb P_p | T_{p-\gamma}[\mathcal L  Q] = \mathbf 0 \}
\text{ and }
\begin{array}{rll}
\mathcal D_p:& \mathbb P_p& \to  \mathbb P_{p-\gamma}\\
& Q&\mapsto T_{p-\gamma}[\mathcal L  Q] 
\end{array}
$$
so $\mathbb Q\mathbb T_p =\ker \mathcal D_p$. In this scalar case, thanks to the graded structure of polynomial spaces, the operator $\mathcal D_p|_{\bigoplus_{\ell=\gamma}^p\widetilde{\mathbb P}_{\ell}}$ was proved to have a block-triangular structure, and these diagonal blocks correspond to the graded operators of degree $-\gamma$ defined by $\displaystyle \sum_{ |\mathbf j|=\gamma} c_{\mathbf j}\left( {\mathbf x_0}  \right) \partial_{{\mathbf x}}^{\mathbf j}$ between spaces of homogeneous polynomials. In the scalar case, such operators were proved to be surjective. From there, a natural block forward substitution procedure lead to an algorithm to construct quasi-Trefftz functions and consequently to the dimension of the quasi-Trefftz spaces.
As a reminder, quasi-Trefftz spaces are smaller than standard polynomial spaces, and the dimension of these discrete spaces with the same order $p$ of best approximation properties in dimension $d$, between the polynomial spaces $\mathbb P_p$ and some quasi-Trefftz spaces $ \mathbb Q\mathbb T_p$ can be described as follows.

{\renewcommand{\arraystretch}{1.75}%
    \begin{tabular}{c|ccc} 
   &   \textbf{$d=2$} & \textbf{$d=3$} & \textbf{$d$} \\
      \hline
   $ \mathbb P_p$&  $\displaystyle\frac{(p+1)(p+2)}{2}=O\left(p^2\right)$ & $\displaystyle\frac{(p+1)(p+2)(p+3)}{6}=O\left(p^3\right)$ & $O(p^{d})$ \\
   $ \mathbb Q\mathbb T_p$ with $\gamma =2$&  $2p+1=O\left(p\right)$ & $(p+1)^2=O\left(p^2\right)$& $O(p^{d-1})$  \\
   $ \mathbb Q\mathbb T_p$ with $\gamma =3$&  $3p=O\left(p\right)$ & $\frac{3p^2+3p+2}2=O\left(p^2\right)$& $O(p^{d-1})$ 
    \end{tabular}}
    
    For systems of PDEs, the challenge comes from the coupling between components within the highest order derivatives in the differential operator. If these components are not coupled, then extending the scalar framework is straightforward, but for Equation \eqref{eq:Max} this is not the case: the curl-curl operator is not surjective between spaces of homogeneous polynomials. 
Yet the quasi-Trefftz operator retains an underlying triangular structure, and it is thanks to the quasi-divergence condition that the quasi-Trefftz spaces will be expressed -as in the scalar case- as the kernel of different operators with an underlying triangular structure and surjective diagonal blocks.
The present work leverages the block structure of this operator for a detailed study of deRham complex operators between spaces of homogeneous polynomials, and introduces a Helmholtz decomposition of homogeneous polynomial vector fields. 
{\bf The outcome of this work is a notion of quasi-Trefftz spaces with same best approximation properties as the corresponding standard polynomial spaces but with an optimal dimension, as well as a procedure to construct a quasi-Trefftz basis.}
 
 \section{
Differential operators and exact sequences
 }
 \label{sec:GDC}

Consider the gradient, divergence and curl operators between standard spaces of polynomials, namely:
$$
\begin{array}{rrcl}
Grad_p:&\mathbb P_{p+1} &\to&(\mathbb P_p)^3,\\
&  \Phi &\mapsto& \nabla \Phi,
\end{array}
\quad 
\begin{array}{rrcl}
Div_p:&(\mathbb P_{p+1})^3 &\to&\mathbb P_p,\\
& \mathbf \Phi &\mapsto& \nabla \cdot\mathbf \Phi,
\end{array}
\text{ and }
\begin{array}{rrcl}
Curl_p:&(\mathbb P_{p+1})^3 &\to&(\mathbb P_{p})^3,\\
& \mathbf \Phi &\mapsto& \nabla \times\mathbf \Phi.
\end{array}
$$
 These maps are graded linear maps of degree $-1$ on polynomials, 
 and for each of these, the range and kernel, their dimensions and a basis of each are of interest to construct a quasi-Trefftz basis for Maxwell's equation. 
Also consider two Laplacian operators between standard spaces of polynomials:
$$
\begin{array}{rrcl}
Lap_p:&\mathbb P_{p+2} &\to&\mathbb P_p,\\
&  \Phi &\mapsto& \Delta \Phi,
\end{array}
\text{ and }
\begin{array}{rrcl}
\overrightarrow{Lap}_p:&(\mathbb P_{p+2} )^3&\to&(\mathbb P_p)^3,\\
&  \mathbf\Phi &\mapsto& \overrightarrow\Delta \mathbf\Phi,
\end{array}
$$
Both of these operators are scalar-valued graded operators of degree -2. 
The former operator, as it is defined between spaces of scalar polynomials, falls under the general framework developed in \cite{imbertgerard2025localtaylorbasedpolynomialquasitrefftz}. The latter is crucial in the study of Equation \eqref{eq:Max} as $\nabla\times\nabla\times \mathbf V = -\overrightarrow \Delta \mathbf V +\nabla(\nabla\cdot  \mathbf V)$.

Each of these operators is of the generic form
\begin{equation}\label{eq:GenOp}
\mathcal T:\bigoplus_{l=0}^{p+\gamma} \left(\widetilde{\mathbb P}_l\right)^n\to\bigoplus_{k=0}^{p} \left(\widetilde{\mathbb P}_k\right)^m ,\text{ for }n,m\in\{1,3\},
\end{equation}
where $\gamma\in\{1,2\}$ is the order of $\mathcal T$, and these graded operators can be decomposed as follows: for all $\left\{x_l\in \left(\widetilde{\mathbb P}_l\right)^n, l\in[\![0,p+\gamma]\!]\right\}$,
\begin{equation}\label{eq:GenOpStr}
\mathcal T\left(\sum_{k=0}^{p+\gamma} x_k\right) = \mathcal T_*\left(\sum_{k=0}^{\gamma-1} x_k\right)+\sum_{k=0}^{p}\mathcal T_{k}\left( x_k\right)
\text{ with }
\mathcal T_*: ({\mathbb P}_{\gamma-1})^n\to \{(\mathbf 0)^m\}
\text{ and }
\mathcal T_{k}:\left(\widetilde{\mathbb P}_{k+\gamma}\right)^n\to \left(\widetilde{\mathbb P}_k\right)^m
\end{equation}
Moreover, since $\mathcal R(\mathcal T_*) = \{\mathbf 0^m\}$, then
\begin{equation}\label{eq:Rs}
\mathcal R(\mathcal T) =  \bigoplus_{0\leq k\leq p} \mathcal R(\mathcal T_{k})
\end{equation}
Hence, the graded structure of polynomial spaces can be leveraged to study the kernels and ranges of the operators $\mathcal T_k$ and $\mathcal T_*$.

This preliminary section focuses on basic properties of these various differential operators that form building blocks for the rest of the article.

\subsection{First order operators 
}
\label{ssec:hom}

For each of the gradient, divergence, and curl operators, the operators $\mathcal T_k$ and $\mathcal T_*$ are defined, and their kernels and ranges are described by a basis.


For the gradient operator, the input is scalar-valued and the output is vector-valued; in the generic form \eqref{eq:GenOp} this corresponds to objects of dimensions $n=1$ for its domain and $m=3$ for its co-domain.
The next definition follows the decomposition \eqref{eq:GenOpStr}.
\begin{dfn}
For any $k\in\mathbb N_0$ and any $\mathbf x_0\in\mathbb R^3$, the gradient operator between spaces of homogeneous polynomials is defined as
\begin{equation}\label{eq:defGk}
\begin{array}{rrcl}
G_k:&\widetilde{\mathbb P}_{k+1} &\to&\left(\widetilde{\mathbb P}_k\right)^3,\\
& \Phi &\mapsto& \nabla \Phi,
\end{array}
\text{ and }
\begin{array}{rrcl}
G_*:&\widetilde{\mathbb P}_{0} &\to&\{\mathbf 0^3\},\\
& \Phi &\mapsto&\mathbf 0^3.
\end{array}
\end{equation}
\end{dfn}
The corresponding kernel and range can be described as follows.
\begin{lmm}\label{lmm:gradk}
For any $k\in\mathbb N_0$ and any $\mathbf x_0\in\mathbb R^3$,
$ \ker (G_k) = \{ \mathbf 0\}$ while the set
$$
\left\{
\begin{pmatrix}
i_1(\mathbf X-\mathbf x_0)^{\mathbf i-\mathbf e_1},
i_2(\mathbf X-\mathbf x_0)^{\mathbf i-\mathbf e_2},
i_3(\mathbf X-\mathbf x_0)^{\mathbf i-\mathbf e_3}
\end{pmatrix}, \mathbf i\in\mathbb N_0^3, |\mathbf i|=k+1
\right\},
$$
forms a basis of $\mathcal R(G_k)$, so the dimensions of the kernel and range of $G_k$ are
$$
\dim\ker(G_k)  = 0
\text{ and }
rk(G_k) = \frac 12(k+2)(k+3).
$$
\end{lmm}

\begin{proof}
Only constant functions have a zero gradient. Moreover, for any $k\in\mathbb N_0$ $ \ker (G_k)\subset \widetilde{\mathbb P}_{k+1}$, hence $ \ker (G_k) = \{ \mathbf 0\}$, and so $\dim\ker(G_k)  = 0$.

Since $\dim \widetilde{\mathbb P}_{k+1} =  \frac 12(k+2)(k+3)$, the rank-nullity theorem gives $rk(G_k) = \frac 12(k+2)(k+3)$. Besides, 
a basis for $\mathcal R(G_k)$ is given by
$$
\left\{
\begin{pmatrix}
i_1(\mathbf X-\mathbf x_0)^{\mathbf i-\mathbf e_1},
i_2(\mathbf X-\mathbf x_0)^{\mathbf i-\mathbf e_2},
i_3(\mathbf X-\mathbf x_0)^{\mathbf i-\mathbf e_3}
\end{pmatrix}, \mathbf i\in\mathbb N_0^3, |\mathbf i|=k+1
\right\},
$$
where for any $g\in\{1,2,3\}$, $0(\mathbf X-\mathbf x_0)^{\mathbf i-\mathbf e_g}=\mathbf 0$ is the zero scalar polynomial.
Indeed, it is a set of $12(k+2)(k+3)$ elements belonging to the range, and they are linearly independent because
$$
\sum_{|\mathbf i|=k+1} \alpha_{\mathbf i} \nabla (\mathbf X-\mathbf x_0)^{\mathbf i} = \mathbf 0^3
\Rightarrow
\nabla\left(\sum_{|\mathbf i|=k+1} \alpha_{\mathbf i} (\mathbf X-\mathbf x_0)^{\mathbf i} \right)= \mathbf 0^3
\Rightarrow
\sum_{|\mathbf i|=k+1} \alpha_{\mathbf i} (\mathbf X-\mathbf x_0)^{\mathbf i} = \mathbf 0,
$$
which in turn implies that $ \alpha_{\mathbf i} = 0$ for all $\mathbf i $ such that $|\mathbf i|=k+1$
since $\{(\mathbf X-\mathbf x_0)^{\mathbf i}, \mathbf i\in\mathbb N_0^3, |\mathbf i|=k+1\}$ forms a basis of $\widetilde{\mathbb P}_{k+1} $.
\end{proof}


For the divergence operator, the input is vector-valued and the output is scalar-valued; in the generic form \eqref{eq:GenOp} this corresponds to objects of dimensions $n=3$ for its domain and $m=1$ for its co-domain.
The next definition follows again the decomposition \eqref{eq:GenOpStr}.
\begin{dfn}
For any $k\in\mathbb N_0$ and any $\mathbf x_0\in\mathbb R^3$, the divergence operator between spaces of homogeneous polynomials is defined as
\begin{equation}\label{eq:defDk}
\begin{array}{rrcl}
D_k:&\left(\widetilde{\mathbb P}_{k+1}\right)^3 &\to&\widetilde{\mathbb P}_k,\\
&\mathbf  \Phi &\mapsto& \nabla \cdot\mathbf \Phi,
\end{array}
\text{ and }
\begin{array}{rrcl}
D_*:&\left(\widetilde{\mathbb P}_{0}\right)^3&\to&\{\mathbf 0\},\\
&\mathbf \Phi &\mapsto&\mathbf 0.
\end{array}.
\end{equation}
\end{dfn}
The corresponding kernel and range can be described as follows.
\begin{lmm}\label{lmm:divk}
For any $k\in\mathbb N_0$ and any $\mathbf x_0\in\mathbb R^3$,
the union of the following sets forms a basis for the kernel $\ker(D_k) $:
\begin{itemize}
\item $\{\begin{pmatrix}(\mathbf X-\mathbf x_0)^{\mathbf i},\mathbf 0,\mathbf 0\end{pmatrix}, \forall \mathbf i\in\mathbb N^3$ such that $|\mathbf i|=k+1,i_1 = 0\}$,
\item $\{\begin{pmatrix}\mathbf 0,(\mathbf X-\mathbf x_0)^{\mathbf i},\mathbf 0\end{pmatrix},\forall\mathbf i\in\mathbb N^3$ such that $|\mathbf i|=k+1,i_2 = 0\}$,
\item $\{\begin{pmatrix}\mathbf 0,\mathbf 0,(\mathbf X-\mathbf x_0)^{\mathbf i}\end{pmatrix},\forall\mathbf i\in\mathbb N^3$ such that $|\mathbf i|=k+1,i_3 = 0\}$,
\item $\{\begin{pmatrix}(i_3+1)(\mathbf X-\mathbf x_0))^{\mathbf i},\mathbf 0,-i_1(\mathbf X-\mathbf x_0)^{\mathbf i-\mathbf e_1 +\mathbf e_3}\end{pmatrix},\forall \mathbf i\in\mathbb N^3$ such that $|\mathbf i|=k+1,i_1> 0\}$,
\item $\{\begin{pmatrix}\mathbf 0,(i_3+1)(\mathbf X-\mathbf x_0)^{\mathbf i},-i_2(\mathbf X-\mathbf x_0)^{\mathbf i-\mathbf e_2 +\mathbf e_3}\end{pmatrix}\forall\mathbf i\in\mathbb N^3$ such that $|\mathbf i|=k+1,i_2> 0\}$,
\end{itemize}
while the operator $D_k$ is surjective,
so the dimensions of the kernel and range of $D_k$ are
$$
\dim\ker(D_k)  = (k+2)(k+4)
\text{ and }
rk(D_k) = \frac 12(k+1)(k+2).
$$
\end{lmm}
\begin{proof}
All elements of the canonical basis of $\widetilde{\mathbb P}_k$, namely $\{(\mathbf X-\mathbf x_0)^{\mathbf i}, \mathbf i\in\mathbb N_0^3, |\mathbf i|=k\}$, belong to the range since
$$\forall \mathbf i\in\mathbb N_0^3, |\mathbf i|=k,
\nabla\cdot\begin{pmatrix}
\frac 1{i_1+1}(\mathbf X-\mathbf x_0)^{\mathbf i+\mathbf e_1},\mathbf 0,\mathbf 0
\end{pmatrix}
=(\mathbf X-\mathbf x_0)^{\mathbf i}
\text{ where }
\begin{pmatrix}
\frac 1{i_1+1}(\mathbf X-\mathbf x_0)^{\mathbf i+\mathbf e_1},\mathbf 0,\mathbf 0
\end{pmatrix}\in\left(\widetilde{\mathbb P}_{k+1}\right)^3.
$$
Hence the operator is surjective and $rk(D_k) =\dim\widetilde{\mathbb P}_k$.

Then from the rank-nullity theorem, 
$$
\dim\ker(D_k) = \frac 32(k+2)(k+3) - \frac 12 (k+1)(k+2) = (k+2)(k+4).
$$
Consider the candidate union of sets. The cardinals of the first three sets  are each equal to $k+2$,
while
the cardinals of the last two sets are both equal to 
$\dim\widetilde{\mathbb P}_{k+1}-(k+2) 
= \frac 12 (k+2)(k+3)-(k+2)
$. 
So together they form a set of $  3(k+2)+2\frac 12(k+2)(k+1) = (k+2)(k+4)$ elements that are clearly linearly independent, hence they form a basis of the space they span. Moreover, this space is included in $\ker(D_k) $ and has the same dimension. This proves the claim. 
\end{proof}
Further useful properties of the divergence operator will be studied in Section \ref{ssec:restr}.


For the curl operator, the input is vector-valued and the output is vector-valued; in the generic form \eqref{eq:GenOp} this corresponds to objects of dimensions $n=3$ for its domain and $m=3$ for its co-domain.
One more, the next definition follows the decomposition \eqref{eq:GenOpStr}.
\begin{dfn} 
For any $k\in\mathbb N_0$ and any $\mathbf x_0\in\mathbb R^3$, the curl operator between spaces of homogeneous polynomials is defined as
\begin{equation}\label{eq:defCk}
\begin{array}{rrcl}
C_k:&\left(\widetilde{\mathbb P}_{k+1}\right)^3 &\to&\left(\widetilde{\mathbb P}_k\right)^3,\\
&\mathbf  \Phi &\mapsto& \nabla \times\mathbf \Phi,
\end{array}
\text{ and }
\begin{array}{rrcl}
C_*:&\left(\widetilde{\mathbb P}_{0}\right)^3&\to&\{\mathbf 0^3\},\\
&\mathbf \Phi &\mapsto&\mathbf 0^3.
\end{array}
\end{equation}
\end{dfn}
The corresponding kernel and range can be described as follows.
\begin{lmm}\label{lmm:curlk}
For any $k\in\mathbb N$ and any $\mathbf x_0\in\mathbb R^3$,
the union of the following sets forms a basis for the range $\mathcal R(C_k) $:
\begin{itemize}
\item $\{\begin{pmatrix}(\mathbf X-\mathbf x_0)^{\mathbf i},\mathbf 0,\mathbf 0\end{pmatrix}, \forall \mathbf i\in\mathbb N^3$ such that $|\mathbf i|=k,i_1 = 0\}$,
\item $\{\begin{pmatrix}\mathbf 0,(\mathbf X-\mathbf x_0)^{\mathbf i},\mathbf 0\end{pmatrix},\forall\mathbf i\in\mathbb N^3$ such that $|\mathbf i|=k,i_2 = 0\}$,
\item $\{\begin{pmatrix}\mathbf 0,\mathbf 0,(\mathbf X-\mathbf x_0)^{\mathbf i}\end{pmatrix},\forall\mathbf i\in\mathbb N^3$ such that $|\mathbf i|=k,i_3 = 0\}$,
\item $\{\begin{pmatrix}(i_3+1)(\mathbf X-\mathbf x_0)^{\mathbf i},\mathbf 0,-i_1(\mathbf X-\mathbf x_0)^{\mathbf i-\mathbf e_1 +\mathbf e_3}\end{pmatrix},\forall \mathbf i\in\mathbb N^3$ such that $|\mathbf i|=k,i_1> 0\}$,
\item $\{\begin{pmatrix}\mathbf 0,(i_3+1)(\mathbf X-\mathbf x_0)^{\mathbf i},-i_2(\mathbf X-\mathbf x_0)^{\mathbf i-\mathbf e_2 +\mathbf e_3}\end{pmatrix}\forall\mathbf i\in\mathbb N^3$ such that $|\mathbf i|=k,i_2> 0\}$;
\end{itemize}
$\{(\mathbf X-\mathbf x_0)^{\mathbf i}, \mathbf i\in\mathbb N_0^3, |\mathbf i|=0\}$ forms a basis of $\mathcal R(C_0) =({\mathbb P}_{0})^3$,
while for any $k\in\mathbb N_0$ a basis of $\ker(C_{k}) $ is given by $\left\{
\begin{pmatrix}
i_1(\mathbf X-\mathbf x_0)^{\mathbf i-\mathbf e_1},
i_2(\mathbf X-\mathbf x_0)^{\mathbf i-\mathbf e_2},
i_3(\mathbf X-\mathbf x_0)^{\mathbf i-\mathbf e_3}
\end{pmatrix}, \mathbf i\in\mathbb N_0^3, |\mathbf i|=k+2
\right\}$.

Hence for any $k\in\mathbb N_0$, the dimensions of the kernel and range of $C_k$ are
$$
\dim\ker(C_k)  = {(k+3)(k+4)}/2
\text{ and }
rk(C_k) = (k+1)(k+3).
$$
\end{lmm}
\begin{proof}
It is a classical results that for any function $\mathbf f\in(\mathcal C^1)^3$, $\nabla\cdot(\nabla\times\mathbf f) =0$, and this in particular implies that $ \mathcal R(C_0)\subset \ker(D_*)$ and $ \mathcal R(C_k) \subset \ker(D_{k-1})$ if $k>0$.
These inclusions are actual identities, indeed, on the one hand,
\begin{itemize}
\item $\forall \mathbf i\in\mathbb N_0^3$ such that $|\mathbf i|=k+1,i_1 = 0$, $\begin{pmatrix}(\mathbf X-\mathbf x_0)^{\mathbf i},\mathbf 0,\mathbf 0\end{pmatrix} = \nabla\times \begin{pmatrix}\mathbf 0,\mathbf 0,\frac1{i_2+1}(\mathbf X-\mathbf x_0)^{\mathbf i+\mathbf e_2}\end{pmatrix}$,
\item $\forall\mathbf i\in\mathbb N_0^3$ such that $|\mathbf i|=k+1,i_2 = 0$, $\begin{pmatrix}\mathbf 0,(\mathbf X-\mathbf x_0)^{\mathbf i},\mathbf 0\end{pmatrix} = \nabla\times \begin{pmatrix}\frac1{i_3+1}(\mathbf X-\mathbf x_0)^{\mathbf i+\mathbf e_3,\mathbf 0,\mathbf 0}\end{pmatrix}$,
\item $\forall\mathbf i\in\mathbb N_0^3$ such that $|\mathbf i|=k+1,i_3 = 0$, $\begin{pmatrix}\mathbf 0,\mathbf 0,(\mathbf X-\mathbf x_0)^{\mathbf i}\end{pmatrix} =  \nabla\times \begin{pmatrix}\mathbf 0,\mathbf 0,\frac1{i_1+1}(\mathbf X-\mathbf x_0)^{\mathbf i+\mathbf e_1}\end{pmatrix}$,
\item $\forall \mathbf i\in\mathbb N_0^3$ such that $|\mathbf i|=k+1,i_1> 0$, $\begin{pmatrix}-(i_3+1)(\mathbf X-\mathbf x_0)^{\mathbf i},\mathbf 0,i_1(\mathbf X-\mathbf x_0)^{\mathbf i-\mathbf e_1 +\mathbf e_3}\end{pmatrix} = \nabla\times \begin{pmatrix}\mathbf 0,(\mathbf X-\mathbf x_0)^{\mathbf i+\mathbf e_3},\mathbf 0\end{pmatrix}$,
\item $\forall\mathbf i\in\mathbb N_0^3$ such that $|\mathbf i|=k+1,i_2> 0$, $\begin{pmatrix}\mathbf 0,(i_3+1)(\mathbf X-\mathbf x_0)^{\mathbf i},-i_2(\mathbf X-\mathbf x_0)^{\mathbf i-\mathbf e_2 +\mathbf e_3}\end{pmatrix} =  \nabla\times \begin{pmatrix}(\mathbf X-\mathbf x_0)^{\mathbf i+\mathbf e_3,\mathbf 0,\mathbf 0}\end{pmatrix}$.
\end{itemize}
So if $k\in\mathbb N$ for any $\mathbf P$ in the basis of $ \ker(D_{k-1})$ provided in Lemma \ref{lmm:divk}, there exists $\mathbf Q\in\left(\widetilde{\mathbb P}_{k+1}\right)^3$ such that $\nabla\times \mathbf Q = \mathbf P$, or equivalently $\mathbf P\in\mathcal R(C_k) $: hence $ \ker(D_{k-1})\subset \mathcal R(C_k) $. On the other hand,
$$\begin{pmatrix}(\mathbf X-\mathbf x_0)^{\mathbf 0},\mathbf 0,\mathbf 0\end{pmatrix} = \nabla\times \begin{pmatrix}\mathbf 0,\mathbf 0,(\mathbf X-\mathbf x_0)^{\mathbf e_2}\end{pmatrix}, \quad
\begin{pmatrix}\mathbf 0,(\mathbf X-\mathbf x_0)^{\mathbf 0},\mathbf 0\end{pmatrix} = \nabla\times \begin{pmatrix}(\mathbf X-\mathbf x_0)^{\mathbf e_3,\mathbf 0,\mathbf 0}\end{pmatrix}, 
$$
$$\text{ and }
\begin{pmatrix}\mathbf 0,\mathbf 0,(\mathbf X-\mathbf x_0)^{\mathbf 0}\end{pmatrix} =  \nabla\times \begin{pmatrix}\mathbf 0,\mathbf 0,(\mathbf X-\mathbf x_0)^{\mathbf e_1}\end{pmatrix},$$
so $\ker(D_*)=({\mathbb P}_{0})^3\subset \mathcal R(C_0)$.
This proves the claims about ranges. 

It is also a classical result that for any function $ f\in\mathcal C^1$, $\nabla\times(\nabla f) =0$, and this in particular implies that $ \mathcal R(G_{k+1}) \subset \ker(C_{k})$. Moreover, from the rank-nullity theorem, 
$$
\dim  \ker(C_{k}) =  \frac 32(k+2)(k+3) - (k+1)(k+3),
$$
hence, according to Lemma \ref{lmm:gradk}, it is equal to $rk(G_{k+1})={(k+3)(k+4)}/2$.
This proves the claims about kernels. 
\end{proof}

\begin{rmk}\label{rmk:ESk}
Together, these operators form the following exact sequence:
$$
\{ \mathbf 0\}
\xrightarrow[]{}\widetilde{\mathbb P}_{k+3}
\xrightarrow[]{G_{k+2}}\left(\widetilde{\mathbb P}_{k+2}\right)^3
\xrightarrow[]{C_{k+1}}\left(\widetilde{\mathbb P}_{k+1}\right)^3
\xrightarrow[]{D_k}\widetilde{\mathbb P}_{k}
\xrightarrow[]{} \{ \mathbf 0\}.
$$
This directly follows from Lemmas \ref{lmm:gradk}, \ref{lmm:divk} and \ref{lmm:curlk} since
$$
\ker (G_{k+2}) = \{\mathbf 0\},\
\mathcal R(G_{k+2}) = \ker(C_{k+1}),\
\mathcal R(C_{k+1}) = \ker(D_{k}) \text{ and }
\mathcal R(D_{k})=\widetilde{\mathbb P}_{k}.
$$
\end{rmk}

\subsection{Second order operators}

For the Laplacian operator, the input is scalar-valued and the output is scalar-valued; in the generic form \eqref{eq:GenOp} this corresponds to objects of dimensions $n=1$ for its domain and $m=1$ for its co-domain.
The next definition follows again the decomposition \eqref{eq:GenOpStr}.
\begin{dfn}
For any $k\in\mathbb N_0$ and any $\mathbf x_0\in\mathbb R^3$, the (scalar) Laplacian operator between spaces of homogeneous polynomials is defined as
\begin{equation*}
\begin{array}{rrcl}
L_k:&\widetilde{\mathbb P}_{k+2} &\to&\widetilde{\mathbb P}_k,\\
&\Phi &\mapsto& \Delta \Phi,
\end{array}
\text{ and }
\begin{array}{rrcl}
L_*:&{\mathbb P}_{1}&\to&\{\mathbf 0\},\\
&\Phi &\mapsto&\mathbf 0.
\end{array}.
\end{equation*}
\end{dfn}
The corresponding kernel and range can be described as follows.
\begin{prop}\label{prop:Lapk}
For any $k\in\mathbb N_0$ and any $\mathbf x_0\in\mathbb R^3$, the operator $L_k$ is surjective,
so the dimensions of the kernel and range of $L_k$ are
$$
\dim\ker(L_k)  = 2k+5 
\text{ and }
rk(L_k) = \frac 12(k+1)(k+2).
$$
\end{prop}
\begin{proof}
As proved in \cite{imbertgerard2025localtaylorbasedpolynomialquasitrefftz},
$L_k|_{\mathbb V_{k+2}}$ is bijective if $\mathbb V_{k+2}:=Span\{(\mathbf X-\mathbf x_0)^{2\mathbf e_1+\mathbf i},\mathbf i\in\mathbb N_0^d,|\mathbf i|=k\}$, and it is clear that $\dim \mathbb V_{k+2} = \dim \widetilde{\mathbb P}_k$.
Therefore 
$\dim rk(L_k) =\dim\widetilde{\mathbb P}_k$, and from the rank-nullity theorem $\dim\ker(L_k)  = \dim\widetilde{\mathbb P}_{k+2}-\dim\widetilde{\mathbb P}_k$.
\end{proof}

For the so-called vector Laplacian operator, both the input and the output are vector-valued, and the output in the component-wise (scalar) Laplacian of the input; in the generic form \eqref{eq:GenOp} this corresponds to objects of dimensions $n=3$ for its domain and $m=3$ for its co-domain.
The next definition follows again the decomposition \eqref{eq:GenOpStr}.
\begin{dfn}
For any $k\in\mathbb N_0$ and any $\mathbf x_0\in\mathbb R^3$, the vector Laplacian operator between spaces of homogeneous polynomials is defined as
\begin{equation*}
\begin{array}{rrcl}
\overrightarrow L_k:&\left(\widetilde{\mathbb P}_{k+2}\right)^3 &\to&\left(\widetilde{\mathbb P}_k\right)^3,\\
&\mathbf  \Phi =(\Phi_1, \Phi_2,\Phi_3)&\mapsto&( \Delta \Phi_1,\Delta \Phi_2,\Delta \Phi_3),
\end{array}
\text{ and }
\begin{array}{rrcl}
\overrightarrow L_*:&\left({\mathbb P}_{1}\right)^3&\to&\{\mathbf 0^3\},\\
&\mathbf \Phi &\mapsto&\mathbf 0^3.
\end{array}
\end{equation*}
\end{dfn}

\begin{rmk}
The operator $\overrightarrow L_k$ has a block upper-triangular structure. Indeed, from bases of $\ker(D_{k+1})$ and $\ker(D_{k-1})$ introduced in Lemma \ref{lmm:divk}, define
$$\mathbb A_1^k:=Span\{\begin{pmatrix}(\mathbf X-\mathbf x_0)^{\mathbf i},\mathbf 0,\mathbf 0\end{pmatrix}, \forall \mathbf i\in\mathbb N_0^3\text{  such that }|\mathbf i|=k,i_1 = 0\},$$
$$\mathbb A_2^k:=Span\{\begin{pmatrix}\mathbf 0,(\mathbf X-\mathbf x_0)^{\mathbf i},\mathbf 0\end{pmatrix},\forall\mathbf i\in\mathbb N_0^3\text{  such that }|\mathbf i|=k,i_2 = 0\},$$
$$\mathbb A_3^k:=Span\{\begin{pmatrix}\mathbf 0,\mathbf 0,(\mathbf X-\mathbf x_0)^{\mathbf i}\end{pmatrix},\forall\mathbf i\in\mathbb N_0^3\text{  such that }|\mathbf i|=k,i_3 = 0\},$$
$$\mathbb A_4^k:=Span\{\begin{pmatrix}(i_3+1)(\mathbf X-\mathbf x_0)^{\mathbf i},\mathbf 0,-i_1(\mathbf X-\mathbf x_0)^{\mathbf i-\mathbf e_1 +\mathbf e_3}\end{pmatrix},\forall \mathbf i\in\mathbb N_0^3\text{  such that }|\mathbf i|=k,i_1> 0\}, 
$$
$$
\text{ and }
\mathbb A_5^k:=Span\{\begin{pmatrix}\mathbf 0,(i_3+1)(\mathbf X-\mathbf x_0)^{\mathbf i},-i_2(\mathbf X-\mathbf x_0)^{\mathbf i-\mathbf e_2 +\mathbf e_3}\end{pmatrix}\forall\mathbf i\in\mathbb N^3\text{  such that }|\mathbf i|=k,i_2> 0\};$$
then

$$
\overrightarrow L_k(\mathbb A^{k+2}_1)\subset \mathbb A_1^k,\
\overrightarrow L_k(\mathbb A^{k+2}_2)\subset \mathbb A^k_2,\
\overrightarrow L_k(\mathbb A^{k+2}_3)\subset \mathbb A^k_3,$$
$$
\overrightarrow L_k(\mathbb A^{k+2}_4)\subset \mathbb A^k_4\cup \mathbb A^k_1\cup\mathbb A^k_3,\
\text{ and }
\overrightarrow L_k(\mathbb A^{k+2}_5)\subset \mathbb A^k_5\cup \mathbb A^k_2\cup\mathbb A^k_3.
$$
\end{rmk}

Further useful properties of the vector Laplacian operator will also be studied in Section \ref{ssec:restr}.
 
 \section{
Helmholtz decomposition
 }\label{sec:HD}
The spaces of polynomial solenoidal fields and polynomial irrotational fields, both described in the previous section, play an important role in the study of 3D polynomial vector fields.
\begin{dfn}
Consider $k\in\mathbb N$ and $\mathbf x_0\in\mathbb R^3$.
Within the space of homogeneous polynomials of degree $k$, the spaces of polynomial solenoidal fields and polynomial irrotational fields are respectively defined as:
$$
\widetilde{\mathbb S}_k:= \ker(D_{k-1})
\text{ and }
\widetilde{\mathbb I}_k:= \ker(C_{k-1}) .
$$
For the record, all constant vector fields are both solenoidal and irrotational.

\end{dfn}

While the previous material is sufficient to prove the existence of a decomposition of vector fields into the sum of a solenoidal component and an irrotational component, see Section \ref{ssec:adecomp}, such a decomposition is not unique as there are non-trivial polynomials that are both solenoidal and irrotational. Thanks to the introduction of these so-called harmonic polynomials, see Section \ref{ssec:HVF}, polynomial vector fields can be uniquely decomposed into a component that is solenoidal but not irrotational, a component that is irrotational but not solenoidal, and a harmonic component, 
see Section \ref{ssec:HD}

\subsection{Existence of a decomposition}\label{ssec:adecomp}
Homogeneous polynomial vector fields can be decomposed into the sum of a solenoidal component and an irrotational component as follows.
\begin{lmm}
Consider any $k\in\mathbb N_0$ and any $\mathbf x_0\in\mathbb R^3$.
Then $$
\forall \mathbf V \in \left(\widetilde{\mathbb P}_{k+1}\right)^3,
\exists (\mathbf F,\mathbf G)\in  \widetilde{ \mathbb S}_{k+1} \times\widetilde{ \mathbb I}_{k+1}
\text{ such that }
\mathbf V = \mathbf F+\mathbf G.
$$
\end{lmm}
\begin{proof}
Consider any polynomial vector field $\mathbf V\in\left(\widetilde{\mathbb P}_{k+1}\right)^3$. Since $L_k$ is surjective from Proposition \ref{prop:Lapk}, there exists $g\in\widetilde{\mathbb P}_{k+2}$ such that $L_k g = D_k \mathbf V$, which is equivalent to $D_{k}(G_{k+1} g-\mathbf V)=0 $. 
Then $\mathbf F:= \mathbf V-G_{k+1} g$ belongs to $\widetilde{ \mathbb S}_{k+1}$ since $ \ker(D_{k}) =  \mathcal R(C_{k+1})$ according to Lemma \ref{lmm:curlk}, $\mathbf G:= G_{k+1} g$ belongs to $\widetilde{ \mathbb I}_{k+1}$, and they satisfy
$\mathbf V = \mathbf F+\mathbf G$.
\end{proof}
Here is a simple proof of non-uniqueness for this decomposition.
Given , any non-trivial $g\in\ker L_k$ and any $(\mathbf F,\mathbf G)\in\ker(D_k)\times \ker(C_k)$, then
$$
\mathbf F+\mathbf G =( \mathbf F +G_{k+1}g)+(\mathbf G-G_{k+1}g),
$$
where the two sides of this identity are two distinct Helmholtz decompositions of a single vector field.

Moreover, assuming that $(\mathbf F_\alpha,\mathbf G_\alpha)\in\ker(D_k)\times \ker(C_k)$ for $\alpha\in\{1,2\}$ with $\mathbf F_1+\mathbf G_1 = \mathbf F_2+\mathbf G_2$, then the difference $\mathbf F_1-\mathbf F_2 = \mathbf G_2-\mathbf G_1$, denoted hereafter $\mathbf V$, is in $\ker(D_k)\cap \ker(C_k)\subset\left( \widetilde{\mathbb P}_{k+1}\right)^3$.
This motivates the introduction of the notion of harmonic vector fields.

\subsection{
Harmonic vector fields
}\label{ssec:HVF}
The next step towards a unique polynomial Helmholtz decomposition is to consider the intersection of spaces of solenoidal and irrotational vector fields within the space of homogeneous polynomial vector fields.
\begin{dfn}
Consider any $\mathbf x_0\in\mathbb R^3$.
Polynomial harmonic vector fields are polynomial vector fields that are both solenoidal and irrotational, and the corresponding vector spaces of homogeneous polynomials are defined as:
$$
\forall k\in\mathbb N,\
\widetilde{\mathbb H}_k :=\widetilde{ \mathbb S}_k\cap\widetilde{\mathbb I_k}.
$$ 
 In other words, 
$\widetilde{\mathbb H}_k$ is the space of homogeneous polynomial harmonic vector fields of degree $k$.
For the record, the space of constant vector fields, ${\mathbb P}_0^3=\widetilde{\mathbb P}_0^3$, is also equal to
$\widetilde{\mathbb H}_0$.
\end{dfn}
\begin{prop}\label{prop:dimHarm}
Consider any $\mathbf x_0\in\mathbb R^3$. The dimensions of spaces of polynomial harmonic vector fields are given by:
$$
\forall k\in\mathbb N_0, \dim \widetilde{\mathbb H}_k = 2k+3
.
$$
\end{prop}
\begin{proof}
The space $\widetilde{\mathbb H}_0 $ is simply $  \left( \widetilde{\mathbb P}_{0}\right)^3 $, so its dimension is 3.
Considering then $k\in\mathbb N$,
%
%
for any polynomial vector field $\mathbf V\in(\widetilde{\mathbb P}_k)^3$, the condition
$\mathbf V\in\widetilde{\mathbb H}_k$ is equivalent to
$$
\left\{\begin{array}{l}
 \mathbf V \in \mathcal R(G_{k}),\\
 \mathbf V \in \ker(D_{k-1}),
\end{array}\right.
\Leftrightarrow
\left\{\begin{array}{l}
\displaystyle \mathbf V=\sum_{|\mathbf i| = k+1} \alpha_{\mathbf i} \nabla (\mathbf X-\mathbf x_0)^{\mathbf i},\\
\nabla \cdot \mathbf V = \mathbf 0,
\end{array}\right.
$$
$$
\Leftrightarrow
\left[ \mathbf V= \nabla  \left(\sum_{|\mathbf i| = k+1} \alpha_{\mathbf i}(\mathbf X-\mathbf x_0)^{\mathbf i}\right) \text{ and }
\displaystyle\Delta  \left(\sum_{|\mathbf i| = k+1} \alpha_{\mathbf i}(\mathbf X-\mathbf x_0)^{\mathbf i}\right)=\mathbf 0
\right].
$$
So $\widetilde{\mathbb H}_k$ is the space of harmonic gradients, i.e. the space of gradients of elements of $\ker(L_{k-1})$.
Moreover, $\dim\ker(L_{k-1})=2k+3$ from Proposition \ref{prop:Lapk}, and for any basis of $\ker(L_{k-1})$ denoted $\{v_l,l\in[\![1,2k+3]\!]\}\subset\widetilde{\mathbb H}_{k+1}$ then $\{\nabla v_l,l\in[\![1,2k+3]\!]\}\subset\widetilde{\mathbb H}_k$ is linearly independent since 
$$
\sum_{l=1}^{2k+3} \alpha_l \nabla v_l=0
\Rightarrow \sum_{l=1}^{2k+3} \alpha_l  v_l\in \ker G_{k+1} = \{\mathbf 0\}
\Rightarrow \alpha_l=0\ \forall l\in [\![1,2k+3]\!].
$$
Hence $ \dim \nabla  \ker(L_{k-1}) = \dim \ker(L_{k-1})$ so $ \dim \widetilde{\mathbb H}_k = \dim \ker(L_{k-1})$, which proves the first claim.
\end{proof}

By definition, the spaces of harmonic vector fields are subspaces of both spaces of solenoidal fields and irrotational fields,
but there is more.

\subsection{A unique decomposition}\label{ssec:HD}
This section focuses on a unique Helmholtz decomposition for homogeneous polynomial vector fields.
\begin{dfn}\label{def:star}
For any $k\in\mathbb N$ and any $\mathbf x_0\in\mathbb R^3$,
$\widetilde{  \mathbb S}_k^*$ and  $\widetilde{  \mathbb I}_k^*$ refer to the complements of $ \widetilde{\mathbb H}_{k}$ respectively in $\widetilde{ \mathbb S}_k$ and  $\widetilde{\mathbb I_k}$:
$$
\widetilde{ \mathbb S}_k = \widetilde{  \mathbb S}_k^*\oplus \widetilde{\mathbb H}_{k},
\text{ and }
\widetilde{ \mathbb I}_k = \widetilde{  \mathbb I}_k^*\oplus \widetilde{\mathbb H}_{k}.
$$
\end{dfn}
The dimensions of these spaces directly follow from the previous definition.
\begin{cor}\label{cor:dim*s}
For any $k\in\mathbb N$ and any $\mathbf x_0\in\mathbb R^3$, $\dim \widetilde{ \mathbb S}^*_k =k(k+2) $ and $\dim \widetilde{ \mathbb I}^*_k = k(k+1)/2$.
\end{cor}
\begin{proof}
By definition \ref{def:star},
$$
\left\{\begin{array}{l}
\dim \widetilde{ \mathbb S}^*_k=\dim \widetilde{ \mathbb S}_k-\dim \widetilde{\mathbb H}_{k},
\\
\dim \widetilde{ \mathbb I}^*_k=\dim \widetilde{ \mathbb I}_k-\dim \widetilde{\mathbb H}_{k}.
\end{array}\right.
$$
Moreover, 
$\dim \widetilde{ \mathbb S}^*_k
=\dim \widetilde{ \mathbb S}_k-\dim \widetilde{\mathbb H}_{k}
$
and
$\dim \widetilde{ \mathbb I}^*_k
=\dim \widetilde{ \mathbb I}_k-\dim \widetilde{\mathbb H}_{k}
$.
The conclusion then follows from
Lemma \ref{lmm:divk},  as $\dim \widetilde{ \mathbb S}_k = (k+1)(k+3)   $,
Lemma \ref{lmm:curlk}, as $\dim \widetilde{ \mathbb I }_k = (k+2)(k+3)/2$, and
 Proposition \ref{prop:dimHarm}, as $\dim \widetilde{ \mathbb H}_k = 2k+3   $.
\end{proof}

 The following proposition states the existence and uniqueness of the desired Helmholtz decomposition. 
\begin{prop}\label{prop:HDk}
Consider $k\in\mathbb N$ and any $\mathbf x_0\in\mathbb R^3$.
Homogeneous polynomial harmonic vector fields can be decomposed as follows:
$$
\forall \mathbf V \in \left(\widetilde{\mathbb P}_{k}\right)^3,
\exists (\mathbf F,\mathbf G,\mathbf H)\in  \widetilde{ \mathbb S}^*_{k} \times\widetilde{ \mathbb I}^*_{k}\times  \widetilde{\mathbb H}_{k}
\text{ such that }
\mathbf V = \mathbf F+\mathbf G+\mathbf H.
$$
or, equivalently, 
$$
 \left(\widetilde{\mathbb P}_{k}\right)^3
=
\widetilde{ \mathbb S}^*_{k} \oplus\widetilde{ \mathbb I}^*_{k}\oplus \widetilde{\mathbb H}_{k} 
.
$$
Moreover,
$\displaystyle
 \left(\widetilde{\mathbb P}_{0}\right)^3= \widetilde{\mathbb H}_{0} $.
\end{prop}
\begin{proof}
According to Section \ref{ssec:hom} and Proposition \ref{prop:dimHarm}, for homogeneous polynomials,
the dimensions \modTC{of the} kernels of the divergence and curl operators are related to the dimension of the space of harmonic vector fields by:
\begin{equation*}
\dim \widetilde{ \mathbb S}_{k} +\dim \widetilde{ \mathbb I}_{k}- \dim \widetilde{\mathbb H}_{k} 
= \dim\left[  \left(\widetilde{\mathbb P}_{k}\right)^3\right].
\end{equation*}
According to Definition \ref{def:star} and Corollary \ref{cor:dim*s},
this concludes the proof.

\end{proof}

\section{Polynomial quasi-Trefftz spaces}\label{sec:QTspace}
This section returns to the PDE of interest in this work, namely Maxwell's equation \eqref{eq:Max} with the differential operator $\mathcal L_M\coloneqq\nabla\times\nabla\times\cdot -\varepsilon \cdot$. 
In this case, in the sense introduced in Definition \ref{def:atprop}, a natural quasi-Trefftz property for a polynomial function $\mathbf \mathbf Q$ of degree $p$ reads
 \begin{equation}\label{eq:qTMax}
 T_{p-2}\left[\mathcal L_M \mathbf \mathbf Q \right]= \mathbf 0^3
 \end{equation}
 $$
 \Leftrightarrow
\forall k\in[\![0,p-2]\!],C_{k}\circ C_{k+1}[\mathbf \mathbf Q_{k+2}] -(\varepsilon \mathbf \mathbf Q)_k =\mathbf 0^3
 \text{ where }
 \mathbf \mathbf Q = \sum_{k=0}^p\mathbf \mathbf Q_k,\ \mathbf \mathbf Q_k\in\widetilde{\mathbb P}_k\ \forall k\in[\![0,p]\!].
$$
This equivalence simply leverages the graded structure of polynomial spaces, yet the two equivalent properties imply a corresponding divergence condition. Indeed, for any $ \mathbf \mathbf Q \in \left({\mathbb P}_{p}\right)^3$
$$
T_{p-2}[\mathcal L_M \mathbf \mathbf Q] = \mathbf 0^3
\Rightarrow
\nabla\cdot  T_{p-2}\left[ \varepsilon\mathbf \mathbf Q \right]=\mathbf 0,
$$
which follows form $ T_{p-2}\left[ \nabla\times\nabla\times\mathbf \mathbf Q \right]=  \nabla\times\nabla\times\mathbf \mathbf Q $ and $\nabla \cdot ( \nabla\times\nabla\times\mathbf \mathbf Q )=0$, or equivalently
$$
\forall k\in[\![0,p-2]\!],C_{k}\circ C_{k+1}[\mathbf \mathbf Q_{k+2}] -(\varepsilon \mathbf \mathbf Q)_k =\mathbf 0
\Rightarrow
\forall k\in[\![1,p-2]\!],D_{k-1}[(\varepsilon \mathbf \mathbf Q)_k] =\mathbf 0,
$$ 
which follows from $\mathcal R(C_k)=\ker D_{k-1}$.

In the scalar case, $ T_{p-2}\circ\mathcal L|_{\bigoplus_{\ell=2}^p\widetilde{\mathbb P}_{\ell}}$  has a block-triangular structure and its diagonal blocks are surjective.
Here, while the operator $ T_{p-2}\circ\mathcal L_M|_{\bigoplus_{\ell=2}^p\widetilde{\mathbb P}_{\ell}}$ also has a block-triangular structure, its diagonal blocks, corresponding to the operators $C_{k}\circ C_{k+1}$, are not surjective.
However, the corresponding divergence condition can be leveraged to define an operator on $\bigoplus_{\ell=2}^p\widetilde{\mathbb P}_{\ell}$ that has both a block-triangular structure and surjective diagonal blocks.
The Helmholtz decomposition introduced in the previous section then provides a powerful tool to derive an algorithm for the construction of quasi-Trefftz functions.

\subsection{Definition, characterization, best approximation properties}
Besides the natural quasi-Trefftz property \eqref{eq:qTMax}, the following quasi-Trefftz spaces definition includes a quasi-divergence condition.
As discussed above, the quasi-divergence condition is a consequence of the first one up to a lower order, but this definition imposes the two conditions at a higher order.
This divergence condition will be key to studying the corresponding quasi-Trefftz spaces.


\begin{dfn}\label{dfn:QTp}
Given $p\in\mathbb N$ such that $p>2$ and $\mathbf x_0\in\mathbb R^3$,
consider Maxwell's equation \eqref{eq:Max} 
with a scalar variable coefficient 
$\varepsilon\in\mathcal C^{p}(\mathbf x_0)$. 
The associated quasi-Trefftz $\mathbb Q\mathbb T_p$ is defined as
\begin{equation}\label{eq:QTspace}
\mathbb Q\mathbb T_p:=\left\{ 
 \mathbf \mathbf Q \in \left({\mathbb P}_{p}\right)^3,
 T_{p-2}\left[ \nabla\times\nabla\times\mathbf \mathbf Q -\varepsilon \mathbf \mathbf Q \right] = \mathbf 0^3,
  T_{p-1}\left[\nabla\cdot\left(\varepsilon\mathbf \mathbf Q  \right) \right] = \mathbf 0
\right\},
\end{equation}
or equivalently, in terms of homogeneous polynomials,
$$
\begin{array}{rl} \displaystyle
\mathbb Q\mathbb T_p:=& \displaystyle\left\{ 
 \sum_{k=0}^p\mathbf Q_k, 
\mathbf Q_k \in \left(\widetilde{\mathbb P}_{k}\right)^3\  \forall k\in[\![0,p]\!],
C_k\circ C_{k+1}[\mathbf Q_{k+2} ]=\sum_{k'=0}^k\varepsilon_{k-k'} \mathbf Q_{k'} \  \forall k\in[\![0,p-2]\!],\right.\\
&\phantom{\sum}\left. \displaystyle
 D_{k-1}\mathbf Q_{k}  =-\varepsilon_0^{-1} \left( \sum_{k'=0}^{k-1} G_{k-k'-1}\varepsilon_{k-k'}\cdot\mathbf Q_{k'}+\mathbf 1_{k>1}\sum_{k'=1}^{k-1}\varepsilon_{k-k'}D_{k'-1} \mathbf Q_{k'}  \right) \  \forall k\in[\![1,p]\!]
\right\},
 \end{array}
$$
\end{dfn}
\begin{rmk}
The equivalence between the two definitions is due to the fact that: 
$$ 
T_{p-1}\left[  \nabla\cdot\left( \varepsilon\mathbf Q  \right) \right]
  = \sum_{k=1}^{p} D_{k-1}[(\varepsilon\mathbf Q)_k]
\text{ and }
 D_{k-1}[(\varepsilon\mathbf Q)_k]
= \sum_{k'=0}^{k-1}(G_{k-k'-1} \varepsilon_{k-k'})\cdot \mathbf Q_{k'} 
+ \sum_{k'=1}^{k} \varepsilon_{k-k'}D_{k'-1} \mathbf Q_{k'} .$$
\end{rmk}
Furthermore, the space can be conveniently characterized in terms of the different components of the Helmholtz decomposition.

\begin{prop}\label{prop:QTcharact}
Given $p\in\mathbb N$ such that $p>2$ and $\mathbf x_0\in\mathbb R^3$,
consider Maxwell's equation \eqref{eq:Max} 
with scalar variable coefficients $\varepsilon\in\mathcal C^{p}(\mathbf x_0)$. The associated quasi-Trefftz spaces can be described as
$$ 
\begin{array}{l}
\displaystyle
\mathbb Q\mathbb T_p\displaystyle
=\left\{ 
 \sum_{k=0}^p\mathbf Q_k,
 \mathbf Q_0\in (\mathbb P_0)^3,
 \mathbf Q_k=\mathbf F_k+\mathbf G_k+\mathbf H_k \  \forall k\in[\![1,p]\!],
  (\mathbf F_k,\mathbf G_k,\mathbf H_k)\in  \widetilde{ \mathbb S}^*_{k} \times\widetilde{ \mathbb I}^*_{k}\times  \widetilde{\mathbb H}_{k} \  \forall k\in[\![1,p]\!]\right.\\
\displaystyle
\left.
\phantom{\mathbb Q\mathbb T_p=\Bigg\{ bla}
\overrightarrow L_{k} \mathbf F_{k+2} =-\sum_{k'=0}^{k}\varepsilon_{k-k'} \mathbf Q_{k'}
 \  \forall k\in[\![0,p-2]\!],\right.\\
\displaystyle
\left.
\phantom{\mathbb Q\mathbb T_p=\Bigg\{ bla}
D_{k} \mathbf G_{k+1}  =-\varepsilon_0^{-1} \left( \sum_{k'=0}^{k} G_{k-k'}\varepsilon_{k+1-k'}\cdot\mathbf Q_{k'} 
+\mathbf 1_{k>0}\sum_{k'=1}^{k}\varepsilon_{k+1-k'}D_{k'-1} \mathbf G_{k'}  \right) \  \forall k\in[\![0,p-1]\!]
\right\}.
\end{array}
$$
\end{prop}
\begin{proof}
This directly follows from the Helmholtz decomposition introduced in Proposition \ref{prop:HDk} since
$$
  \forall k\in[\![0,p-2]\!],
\forall   (\mathbf F_{k+2},\mathbf G_{k+2},\mathbf H_{k+2})\in  \widetilde{ \mathbb S}^*_{{k+2}} \times\widetilde{ \mathbb I}^*_{{k+2}}\times  \widetilde{\mathbb H}_{{k+2}},
 C_k\circ C_{k+1}[\mathbf F_{k+2}+\mathbf G_{k+2}+\mathbf H_{k+2}]=\overrightarrow L_{k}\mathbf F_{k+2} ,$$
  together with the fact that 
  $$
  \forall k\in[\![0,p-1]\!],\ 
  \forall  (\mathbf F_{k},\mathbf G_{k},\mathbf H_{k})\in  \widetilde{ \mathbb S}^*_{{k}} \times\widetilde{ \mathbb I}^*_{{k}}\times  \widetilde{\mathbb H}_{{k}},\
  D_{k-1}\left[\mathbf F_{k}+\mathbf G_{k}+\mathbf H_{k} \right]
=  D_{k-1}\mathbf G_{k}.
  $$

\end{proof}

This characterization of the quasi-Trefftz spaces 
provides foundations for theoretical developments, including the proof of best approximation properties and the construction of quasi-Trefftz functions thanks to
right inverse operators for the vector Laplacian $\overrightarrow L_k$ and the divergence $D_{k}$ acting on appropriate subspaces of vector-valued polynomials.

\begin{thm}\label{thn:BAP}
Given $p\in\mathbb N$ such that $p>2$ and $\mathbf x_0\in\mathbb R^3$,
consider Maxwell's equation \eqref{eq:Max} 
with scalar variable coefficients $\varepsilon\in\mathcal C^{p}(\mathbf x_0)$, {and} the associated quasi-Trefftz spaces $\mathbb Q\mathbb T_p$. 
{If $\mathbf V \in (\mathcal C^{p+1}(\mathbf x_0))^3$ satisfies the Maxwell's equation \eqref{eq:Max}, then $T_p \left[\mathbf V\right] \in \mathbb Q\mathbb T_p$, and 
\begin{equation*}
    \inf_{\mathbf Q \in \mathbb Q\mathbb T_p} \lVert \mathbf V - \mathbf Q \rVert_{\infty, \mathcal B(\mathbf{x}_0, r)} \leq \frac{6^{p+1}}{(p+1)!} r^{p+1} 
    \lvert \mathbf V \rvert_{\mathcal C^{p+1}( \mathcal B(\mathbf{x}_0, r))}
\end{equation*}
for all $r>0$, where $B(\mathbf{x}_0, r)$ is the ball centered in $\mathbf x_0$ and radius $r$.}
\end{thm}


\begin{proof}
    To show that $T_p \left[\mathbf V\right] \in \mathbb Q\mathbb T_p$, it is sufficient to verify that the Taylor polynomial $T_p \left[\mathbf V\right]$ satisfies the two properties defining the quasi-Trefftz spaces \eqref{eq:QTspace}. Using the properties of Taylor polynomials, we have
    \begin{itemize}
        \item for all multi-index $\mathbf i\in (\mathbb N_0)^3$ such that $\lvert \mathbf i\rvert\leq p-2$
        \begin{equation*}
            \partial^{\mathbf i} \left(\nabla\times\nabla\times T_{p}\left[\mathbf V\right] - \varepsilon T_{p}\left[\mathbf V\right] \right) (\mathbf x_0) = \partial^{\mathbf i} \left(\nabla\times\nabla\times \mathbf V - \varepsilon \mathbf V\right) (\mathbf x_0) = \mathbf 0
        \end{equation*}
        since $\partial^{\mathbf i} \left(\nabla\times\nabla\times T_{p}\left[\mathbf V\right]\right) (\mathbf x_0) = \partial^{\mathbf i} \left(\nabla\times\nabla\times \mathbf V \right) (\mathbf x_0)$, $\partial^{\mathbf i} \left(\varepsilon T_{p}\left[\mathbf V\right]\right) (\mathbf x_0) = \partial^{\mathbf i} \left(\varepsilon V\right) (\mathbf x_0)$, and $\mathbf V$ satisfies the Maxwell's equation \eqref{eq:Max};
        \item for all multi-index $\mathbf i\in (\mathbb N_0)^3$ such that $\lvert \mathbf i\rvert\leq p-1$
        \begin{equation*}
            \partial^{\mathbf i} \nabla\cdot\left(\varepsilon T_{p}\left[\mathbf V\right] \right) (\mathbf x_0) = \partial^{\mathbf i} \nabla\cdot\left(\varepsilon\mathbf V \right) (\mathbf x_0) = \mathbf 0
        \end{equation*}
        since $\mathbf V$ satisfies the Maxwell's equation \eqref{eq:Max}.
    \end{itemize}
    To conclude the proof
    \begin{equation*}
        \begin{split}
            \inf_{\mathbf  Q \in \mathbb Q\mathbb T_p} \lVert \mathbf V - \mathbf Q \rVert_{\infty, \mathcal B(\mathbf{x}_0, r)}
            \leq&\ \lVert \mathbf V - T_{p}\left[\mathbf V\right] \rVert_{\infty, \mathcal B(\mathbf{x}_0, r)}
            \leq \sum_{\substack{\mathbf j \in (\mathbb N_0)^3,\\ \lvert\mathbf j\rvert = p+1}} \frac{1}{\mathbf j!} \sup_{\mathbf x, x_*\in \mathcal B(\mathbf{x}_0, r)} \lVert \partial^{\mathbf j} \mathbf V(\mathbf x_*) (\mathbf x-\mathbf x_0)^{\mathbf j} \rVert \\
            \leq&\ \frac{6^{p+1}}{(p+1)!} r^{p+1} 
            \lvert \mathbf V \rvert_{\mathcal C^{p+1}( \mathcal B(\mathbf{x}_0, r))}.
        \end{split}
    \end{equation*}
\end{proof}


 \subsection{Restricted operators}\label{ssec:restr}
In Proposition \ref{prop:QTcharact}, the operators coming into play are $D_k|_{\tilde{\mathbb I}^*_{k+1}}$ and $\overrightarrow L_k|_{\tilde{\mathbb S}^*_{k+2}}$. This section focuses on proving that both are surjective.
For the divergence operator, it simply follows from the fact that the kernel $ \widetilde{ \mathbb S}_{k+1}$ of the divergence $D_k$.
As for the Laplacian operator, 
$ \widetilde{ \mathbb I}_{k+2}$ is not the kernel of $\overrightarrow L_k$,
yet the result can be proved by considering first the operator $\overrightarrow L_k|_{\tilde{\mathbb S}_{k+2}}$.

\begin{prop}\label{prop:DivRestrict}
Given $\mathbf x_0\in\mathbb R^3$, the restricted divergence operator defined by
$$
\forall k\in\mathbb N_0, \begin{array}{rccc}
D_k|_{\tilde{\mathbb I}^*_{k+1}}:&\tilde{\mathbb I}^*_{k+1}&\to&\widetilde{\mathbb P}_k,\\
&\mathbf  \Phi &\mapsto&\nabla \cdot \mathbf  \Phi ,
\end{array}
$$
 is bijective.
\end{prop}
\begin{proof}
For any $k\in\mathbb N_0$, the operator $D_k$ defined on the domain $\left(\widetilde{\mathbb P}_{k+1}\right)^3$ is surjective according to Lemma \ref{lmm:divk}.
Moreover, according to Proposition \ref{prop:HDk},
$ \left(\widetilde{\mathbb P}_{k+1}\right)^3 = \widetilde{ \mathbb S}_{k+1} \oplus\widetilde{ \mathbb I}^*_{k+1}$.
Since, by definition, $ \widetilde{ \mathbb S}_{k+1}$ is precisely the kernel of $D_k$, then $D_k|_{\tilde{\mathbb I}^*_{k+1}}$ is also surjective.
Besides, from Corollary \ref{cor:dim*s}, $\dim \widetilde{ \mathbb I}^*_{k+1} = (k+1)(k+2)/2$, therefore $\dim \widetilde{ \mathbb I}^*_{k+1} =\dim \widetilde{\mathbb P}_k$. This concludes the proof.

\end{proof}
\begin{rmk}
Since $\tilde{\mathbb I}_{k+1}=\ker ( C_k) = \mathcal R(G_{k+1})$, given $f\in \widetilde{\mathbb P}_k$ then
$$
	\exists \mathbf  \Phi \in\ker(C_{k})=\mathcal R(G_{k+1}), D_k\mathbf  \Phi = f
	\Leftrightarrow
	\exists  \Phi \in\widetilde{\mathbb P}_{k+2}, D_k \circ G_{k+1}  \Phi = f.
$$
Therefore the surjectivity of the restricted divergence operator $D_k|_{\tilde{\mathbb I}_{k+1}}$ directly follows from the surjectivity of the Laplacian operator $L_k$ stated in Proposition \ref{prop:Lapk}.
\end{rmk}

The vector Laplacian acting between spaces of divergence-free fields has the following property.

\begin{prop}\label{prop:VecLapRestrict}
Given $\mathbf x_0\in\mathbb R^3$, the vector Laplacian operator restricted to solenoidal fields, namely
$$
\forall k\in\mathbb N, \begin{array}{rccc}
\overrightarrow L_k|_{\tilde{\mathbb S}_{k+2}}:&\tilde{\mathbb S}_{k+2}&\to&\tilde{\mathbb S}_k,\\
&\mathbf  \Phi 
&\mapsto& \overrightarrow L_k (\mathbf  \Phi )
\end{array}
\text{ or }
\begin{array}{rrcc}
\overrightarrow L_0|_{\tilde{\mathbb S}_2}:&\tilde{\mathbb S}_2&\to&\widetilde{\mathbb P}_{0}^3,\\
&\mathbf  \Phi &\mapsto&\overrightarrow L_0\mathbf  \Phi ,
\end{array}
$$
 is surjective.
As a consequence, 
$$
\forall k\in\mathbb N_0, \dim\ker \left(\overrightarrow L_k|_{\tilde{\mathbb S}_{k+2}}\right) = 4(k+3).
$$
\end{prop}

For the sake of compactness, consider the following notation for elements of divergence-free bases, for any $(k,\mathbf i)\in\mathbb N_0\times \mathbb N_0^3$ with $|\mathbf i|=k$:
\begin{itemize}
\item $\mathbf \Psi^{k,1,\mathbf i}:=\begin{pmatrix}(\mathbf X-\mathbf x_0)^{\mathbf i},\mathbf 0,\mathbf 0\end{pmatrix}$ if $i_1 = 0$;
$\mathbf \Psi^{k,2,\mathbf i}=\begin{pmatrix}\mathbf 0,(\mathbf X-\mathbf x_0)^{\mathbf i},\mathbf 0\end{pmatrix}$ if $i_2 = 0$;\\
$\mathbf \Psi^{k,3,\mathbf i}=\begin{pmatrix}\mathbf 0,\mathbf 0,(\mathbf X-\mathbf x_0)^{\mathbf i}\end{pmatrix}$ if $i_3 = 0$;
\item $\mathbf \Psi^{k,4,\mathbf i}=\begin{pmatrix}(i_3+1)(\mathbf X-\mathbf x_0)^{\mathbf i},\mathbf 0,-i_1(\mathbf X-\mathbf x_0)^{\mathbf i-\mathbf e_1 +\mathbf e_3}\end{pmatrix}$ if $i_1> 0$;\\
$\mathbf \Psi^{k,5,\mathbf i}=\begin{pmatrix}\mathbf 0,(i_3+1)(\mathbf X-\mathbf x_0)^{\mathbf i},-i_2(\mathbf X-\mathbf x_0)^{\mathbf i-\mathbf e_2 +\mathbf e_3}\end{pmatrix}$ if $i_2> 0$.
\end{itemize}
Even though given $\mathbf i\in\mathbb N_0^3$ the value of $k=|\mathbf i|$ is fixed, and therefore the $k$ in this notation is somewhat redundant, this notation will enhance readability in what follows.
\begin{proof}

For $k=0$, it is easy to verify that 
$$
\overrightarrow L_0 \left(\mathbf \Psi^{2,1,2\mathbf e_2}\right)=\mathbf \Psi^{0,1,\mathbf 0}
,\
\overrightarrow L_0 \left(\mathbf \Psi^{2,2,2\mathbf e_3}\right)=\mathbf \Psi^{0,2,\mathbf 0}
\text{ and }
\overrightarrow L_0 \left(\mathbf \Psi^{2,3,2\mathbf e_1}\right)=\mathbf \Psi^{0,3,\mathbf 0}
,
$$
and 
$\{\mathbf \Psi^{0,1,\mathbf 0},\mathbf \Psi^{0,2,\mathbf 0},\mathbf \Psi^{0,3,\mathbf 0}\}$ 
forms a basis of $\widetilde{\mathbb P}_{0}^3$.
Hence $\widetilde{\mathbb P}_{0}^3\subset \overrightarrow L_0\left(\tilde{\mathbb S}_2\right)$.
So the two spaces are equal and the operator is surjective. 

For $k=1$,$\overrightarrow L_1 \left(\tilde{\mathbb S}_3\right)\subset \tilde{\mathbb S}_1$ since $D_{0}\circ \overrightarrow L_1 = L_{0} \circ D_{2}$.
Moreover, 
$$
\overrightarrow L_1\left(  \mathbf \Psi^{3,1,3\mathbf e_3} \right) = 6 \mathbf \Psi^{1,1,\mathbf e_3} ,
\overrightarrow L_1\left(  \mathbf \Psi^{3,1,3\mathbf e_2} \right) = 6 \mathbf \Psi^{1,1,\mathbf e_2} ,
$$
$$
\overrightarrow L_1\left(  \mathbf \Psi^{3,2,3\mathbf e_1} \right) = 6 \mathbf \Psi^{1,2,\mathbf e_1} ,
\overrightarrow L_1\left(  \mathbf \Psi^{3,2,3\mathbf e_3} \right) = 6\mathbf \Psi^{1,2,\mathbf e_3} ,
$$
$$
\overrightarrow L_1\left(  \mathbf \Psi^{3,3,3\mathbf e_1} \right) = 6 \mathbf \Psi^{1,3,\mathbf e_1} ,
\overrightarrow L_1\left(  \mathbf \Psi^{3,3,3\mathbf e_2} \right) = 6 \mathbf \Psi^{1,3,\mathbf e_2} ,
$$
$$
\overrightarrow L_1\left(  \mathbf \Psi^{3,4,3\mathbf e_1} \right) = 6 \mathbf \Psi^{1,4,\mathbf e_1} ,
\overrightarrow L_1\left(  \mathbf \Psi^{3,5,3\mathbf e_2} \right) = 6 \mathbf \Psi^{1,5,\mathbf e_2}, 
$$
and $\{
\mathbf \Psi^{1,1,\mathbf e_3} ,\mathbf \Psi^{1,1,\mathbf e_2} , 
\mathbf \Psi^{1,2,\mathbf e_1} ,\mathbf \Psi^{1,2,\mathbf e_3} ,
 \mathbf \Psi^{1,3,\mathbf e_1} , \mathbf \Psi^{1,3,\mathbf e_2} ,
 \mathbf \Psi^{1,4,\mathbf e_1} ,  \mathbf \Psi^{1,5,\mathbf e_2}
 \}$ forms a basis of $ \tilde{\mathbb S}_1$. Hence $  \tilde{\mathbb S}_1\subset\overrightarrow L_k  \left(\tilde{\mathbb S}_3\right)$. So again the two spaces are equal and the operator is surjective.

For any $k>1$, $\overrightarrow L_k \left( \tilde{\mathbb S}_{k+2}\right)\subset  \tilde{\mathbb S}_{k}$ since $D_{k-1}\circ \overrightarrow L_k = L_{k-1} \circ D_{k+1}$. Furthermore, according to Lemma \ref{lmm:divk} $\dim \tilde{\mathbb S}_{k}=(k+1)(k+3)$, so
$\dim \overrightarrow L_k \left( \tilde{\mathbb S}_{k+2}\right)\geq (k+1)(k+3)$.
Indeed,
 the union of the following sets of polynomial functions in $ \overrightarrow L_k \left( \tilde{\mathbb S}_{k+2}\right)$ is linearly independent:
\begin{itemize}
\item $\{\overrightarrow L_k (\mathbf \Psi^{k+2,1,\mathbf i}), \forall \mathbf i\in\mathbb N_0^3\text{ such that }|\mathbf i|=k+2,i_1 = 0, i_2\leq k\}$ 
with $k+1$ elements,
\item $\{\overrightarrow L_k (\mathbf \Psi^{k+2,2,\mathbf i}), \forall \mathbf i\in\mathbb N_0^3\text{ such that }|\mathbf i|=k+2,i_2 = 0, i_3\leq k\}$
with $k+1$ elements,
\item $\{\overrightarrow L_k (\mathbf \Psi^{k+2,3,\mathbf i}), \forall \mathbf i\in\mathbb N_0^3\text{ such that }|\mathbf i|=k+2,i_3 = 0, i_1\leq k\}$
with $k+1$ elements,
\item $\{\overrightarrow L_k (\mathbf \Psi^{k+2,4,\mathbf i}), \forall \mathbf i\in\mathbb N_0^3\text{ such that }|\mathbf i|=k+2,i_1>2\}$
with $k(k+1)/2$ elements,
\item $\{\overrightarrow L_k (\mathbf \Psi^{k+2,5,\mathbf i}), \forall \mathbf i\in\mathbb N_0^3\text{ such that }|\mathbf i|=k+2,i_2>2\}$
with $k(k+1)/2$ elements,
\end{itemize}
as one can easily verify that 
\begin{itemize}
\item 
$ \overrightarrow L_k ( \mathbf \Psi^{k+2,1,(k+2)\mathbf e_3} )=(k+2)(k+1) \mathbf \Psi^{k,1,k\mathbf e_3}$, 
$ \overrightarrow L_k ( \mathbf \Psi^{k+2,1,\mathbf e_2+(k+1)\mathbf e_3} )=(k+1)k \mathbf \Psi^{k,1,\mathbf e_2+(k-1)\mathbf e_3}$,
and for all $i_2\in[\![2,k]\!]$ $\overrightarrow L_k (\mathbf \Psi^{k+2,1,\mathbf i}) = i_2(i_2-1) \mathbf \Psi^{k+2,1,\mathbf i-2\mathbf e_2}+ i_3(i_3-1)\mathbf \Psi^{k+2,1,\mathbf i-2\mathbf e_3}$,
and together they form a linearly independent set since for any $i_2\in[\![0,k]\!]$ $ \overrightarrow L_k ( \mathbf \Psi^{k+2,1,\mathbf i} )$ is a linear combination of one or two $ \mathbf \Psi^{k,1,\mathbf j} $ with $j_2\leq i_2$;
\item 
$ \overrightarrow L_k ( \mathbf \Psi^{k+2,2,(k+2)\mathbf e_1} )=(k+2)(k+1) \mathbf \Psi^{k,2,k\mathbf e_1}$, 
$ \overrightarrow L_k ( \mathbf \Psi^{k+2,2,\mathbf e_3+(k+1)\mathbf e_1} )=(k+1)k \mathbf \Psi^{k,2,\mathbf e_3+(k-1)\mathbf e_1}$,
and for all $i_3\in[\![2,k]\!]$ $\overrightarrow L_k (\mathbf \Psi^{k+2,2,\mathbf i}) = i_3(i_3-1) \mathbf \Psi^{k+2,2,\mathbf i-2\mathbf e_3}+ i_1(i_1-1)\mathbf \Psi^{k+2,2,\mathbf i-2\mathbf e_1}$,
and together they form a linearly independent set since for any $i_3\in[\![0,k]\!]$ $ \overrightarrow L_k ( \mathbf \Psi^{k+2,2,\mathbf i} )$ is a linear combination of one or two $ \mathbf \Psi^{k,2,\mathbf j} $ with $j_3\leq i_3$;
\item 
$ \overrightarrow L_k ( \mathbf \Psi^{k+2,3,(k+2)\mathbf e_2} )=(k+2)(k+1) \mathbf \Psi^{k,3,k\mathbf e_2}$, 
$ \overrightarrow L_k ( \mathbf \Psi^{k+2,3,\mathbf e_1+(k+1)\mathbf e_2} )=(k+1)k \mathbf \Psi^{k,3,\mathbf e_1+(k-1)\mathbf e_12}$,
and for all $i_1\in[\![2,k]\!]$ $\overrightarrow L_k (\mathbf \Psi^{k+2,3,\mathbf i}) = i_1(i_1-1) \mathbf \Psi^{k+2,3,\mathbf i-2\mathbf e_1}+ i_2(i_2-1)\mathbf \Psi^{k+2,3,\mathbf i-2\mathbf e_2}$,
and together they form a linearly independent set since for any $i_1\in[\![0,k]\!]$ $ \overrightarrow L_k ( \mathbf \Psi^{k+2,3,\mathbf i} )$ is a linear combination of one or two $ \mathbf \Psi^{k,3,\mathbf j} $ with $j_1\leq i_1$,
\item
$ \overrightarrow L_k ( \mathbf \Psi^{k+2,4,(k+2)\mathbf e_1} )=(k+2)(k+1) \mathbf \Psi^{k,4,k\mathbf e_1}$, 
$ \overrightarrow L_k ( \mathbf \Psi^{k+2,4,(k+1)\mathbf e_1+\mathbf e_2} )=(k+1)k \mathbf \Psi^{k,4,(k-1)\mathbf e_1+\mathbf e_2}$, 
 $$\overrightarrow L_k ( \mathbf \Psi^{k+2,4,(k+1)\mathbf e_1+\mathbf e_3} ) =  (k+1)k\mathbf \Psi^{k,4,(k-1)\mathbf e_1+\mathbf e_3}-2(k+1)\mathbf \Psi^{k,3,k\mathbf e_1},$$
 $$\overrightarrow L_k ( \mathbf \Psi^{k+2,4,k\mathbf e_1+\mathbf e_2+\mathbf e_3} ) =  k(k-1)\mathbf \Psi^{k,4,(k-2)\mathbf e_1+\mathbf e_2+\mathbf e_3}-2k\mathbf \Psi^{k,3,(k-1)\mathbf e_1},$$
for all $i_1\in[\![3,k]\!]$ (with $i_2=0$)
  $$\overrightarrow L_k ( \mathbf \Psi^{k+2,4,i_1\mathbf e_1+(k+2-i_1)\mathbf e_3} ) =  i_1(i_1-1)\mathbf \Psi^{k,4,(i_1-2)\mathbf e_1+(k+2-i_1)\mathbf e_3}+(k+3-i_1)(k+2-i_1) \mathbf \Psi^{k,4,i_1\mathbf e_1+(k-i_1)\mathbf e_3}, $$
for all $i_1\in[\![3,k]\!]$ (with $i_3=0$)
  $$\overrightarrow L_k ( \mathbf \Psi^{k+2,4,i_1\mathbf e_1+(k+2-i_1)\mathbf e_2} ) = i_1(i_1-1)\mathbf \Psi^{k,4,(i_1-2)\mathbf e_1+(k+2-i_1)\mathbf e_2}+(k+2-i_1)(k+1-i_1) \mathbf \Psi^{k,4,i_1\mathbf e_1+(k-i_1)\mathbf e_2} ,$$
for all $i_1\in[\![3,k-1]\!]$ (with $i_2=1$)
  $$\overrightarrow L_k ( \mathbf \Psi^{k+2,4,i_1\mathbf e_1+\mathbf e_2+(k+1-i_1)\mathbf e_3} ) 
     =  i_1(i_1-1)\mathbf \Psi^{k,4,(i_1-2)\mathbf e_1+\mathbf e_2+(k+1-i_1)\mathbf e_3} 
     $$$$
    +(k+2-i_1)(k+1-i_1)\mathbf \Psi^{k,4,i_1\mathbf e_1+\mathbf e_2+(k-1-i_1)\mathbf e_3} 
  $$
for all $i_1\in[\![3,k-1]\!]$ (with $i_3=1$)
  $$\overrightarrow L_k ( \mathbf \Psi^{k+2,4,i_1\mathbf e_1+(k+1-i_1)\mathbf e_2+\mathbf e_3} )
   =  i_1(i_1-1)\mathbf \Psi^{k,4,(i_1-2)\mathbf e_1+(k+1-i_1)\mathbf e_2+\mathbf e_3} 
  $$$$
  +(k+1-i_1)(k-i_1)\mathbf \Psi^{k,4,i_1\mathbf e_1+(k-1-i_1)\mathbf e_2+\mathbf e_3} 
-2i_1\mathbf \Psi^{k,3,(i_1-1)\mathbf e_1+(k+1-i_1)\mathbf e_2} ,
  $$
  and finally for all $\mathbf i$ with $i_1>2$, $i_2>1$ and $i_3>1$
  $$
  \overrightarrow L_k ( \mathbf \Psi^{k+2,4,\mathbf i} )
  =i_1(i_1-1) \mathbf \Psi^{k,4,\mathbf i-2\mathbf e_1}
  =i_2(i_2-1) \mathbf \Psi^{k,4,\mathbf i-2\mathbf e_2}
  =i_3(i_3-1) \mathbf \Psi^{k,4,\mathbf i-2\mathbf e_3}, 
  $$
  and together they form a linearly independent set that can be evidenced thanks to the fact that each $ \overrightarrow L_k ( \mathbf \Psi^{k+2,4,\mathbf i} )$ contains a term $ \mathbf \Psi^{k,4,\mathbf i-2\mathbf e_1}$,
\item
$ \overrightarrow L_k ( \mathbf \Psi^{k+2,5,(k+2)\mathbf e_2} )=(k+2)(k+1) \mathbf \Psi^{k,5,k\mathbf e_2}$, 
$ \overrightarrow L_k ( \mathbf \Psi^{k+2,5,(k+1)\mathbf e_2+\mathbf e_1} )=(k+1)k \mathbf \Psi^{k,5,(k-1)\mathbf e_2+\mathbf e_1}$, 
 $$\overrightarrow L_k ( \mathbf \Psi^{k+2,5,(k+1)\mathbf e_2+\mathbf e_3} ) =  (k+1)k\mathbf \Psi^{k,5,(k-1)\mathbf e_2+\mathbf e_3}-2(k+1)\mathbf \Psi^{k,3,k\mathbf e_2},$$
 $$\overrightarrow L_k ( \mathbf \Psi^{k+2,5,\mathbf e_1+k\mathbf e_2+\mathbf e_3} ) =  k(k-1)\mathbf \Psi^{k,5,\mathbf e_1+(k-2)\mathbf e_2+\mathbf e_3}-2k\mathbf \Psi^{k,3,(k-1)\mathbf e_2},$$
for all $i_2\in[\![3,k]\!]$ (with $i_1=0$)
  $$\overrightarrow L_k ( \mathbf \Psi^{k+2,5,i_2\mathbf e_2+(k+2-i_2)\mathbf e_3} ) =  i_2(i_2-1)\mathbf \Psi^{k,5,(i_2-2)\mathbf e_2+(k+2-i_2)\mathbf e_3}+(k+3-i_2)(k+2-i_2) \mathbf \Psi^{k,5,i_2\mathbf e_2+(k-i_2)\mathbf e_3}, $$
for all $i_2\in[\![3,k]\!]$ (with $i_3=0$)
  $$\overrightarrow L_k ( \mathbf \Psi^{k+2,5,(k+2-i_2)\mathbf e_1+i_2\mathbf e_2} ) = i_2(i_2-1)\mathbf \Psi^{k,5,(k+2-i_2)\mathbf e_1+(i_2-2)\mathbf e_2}+(k+2-i_2)(k+1-i_2) \mathbf \Psi^{k,5,(k-i_2)\mathbf e_1+i_2\mathbf e_2} ,$$
for all $i_2\in[\![3,k-1]\!]$ (with $i_1=1$)
  $$\overrightarrow L_k ( \mathbf \Psi^{k+2,5,\mathbf e_1+i_2\mathbf e_2+(k+1-i_2)\mathbf e_3} ) 
     =  i_2(i_2-1)\mathbf \Psi^{k,5,\mathbf e_1+(i_2-2)\mathbf e_2+(k+1-i_2)\mathbf e_3} 
     $$$$
    +(k+2-i_2)(k+1-i_2)\mathbf \Psi^{k,5,\mathbf e_1+i_2\mathbf e_2+(k-1-i_2)\mathbf e_3} 
  $$
for all $i_2\in[\![3,k-1]\!]$ (with $i_3=1$)
  $$\overrightarrow L_k ( \mathbf \Psi^{k+2,5,(k+1-i_2)\mathbf e_1+i_2\mathbf e_2+\mathbf e_3} )
   =  i_2(i_2-1)\mathbf \Psi^{k,5,(k+1-i_2)\mathbf e_1+(i_2-2)\mathbf e_2+\mathbf e_3} 
  $$$$
  +(k+1-i_2)(k-i_2)\mathbf \Psi^{k,5,(k-1-i_2)\mathbf e_1+i_2\mathbf e_2+\mathbf e_3} 
-2i_2\mathbf \Psi^{k,3,(k+1-i_2)\mathbf e_1+(i_2-1)\mathbf e_2} ,
  $$
  and finally for all $\mathbf i$ with $i_1>1$, $i_2>2$ and $i_3>1$
  $$
  \overrightarrow L_k ( \mathbf \Psi^{k+2,5,\mathbf i} )
  =i_1(i_1-1) \mathbf \Psi^{k,5,\mathbf i-2\mathbf e_1}
  =i_2(i_2-1) \mathbf \Psi^{k,5,\mathbf i-2\mathbf e_2}
  =i_3(i_3-1) \mathbf \Psi^{k,5,\mathbf i-2\mathbf e_3}, 
  $$  
  and together they form a linearly independent set that can be evidenced thanks to the fact that each $ \overrightarrow L_k ( \mathbf \Psi^{k+2,5,\mathbf i} )$ contains a term $ \mathbf \Psi^{k,5,\mathbf i-2\mathbf e_2}$
\end{itemize}
This proves the claim that
$\dim \overrightarrow L_k \left(\tilde{\mathbb S}_{k+2}\right)\geq (k+1)(k+3)$
as the union of the previous sets contains $3(k+1)+ 2\frac{k(k+1)}2 = (k+1)(k+3)$ elements.
As a consequence $\overrightarrow L_k \left(\tilde{\mathbb S}_{k+2}\right)=\tilde{\mathbb S}_{k}$ and the operator is surjective.

As a summary, for any $k\in\mathbb N_0$, $\dim\overrightarrow L_k\left({\tilde{\mathbb S}_{k+2}}\right) = (k+1)(k+3)$.
The final result follows from the rank-nullity theorem combined with Lemma \ref{lmm:divk} since 
$$
\forall k\in\mathbb N,
\dim  \tilde{\mathbb S}_{k+2} - \dim \tilde{\mathbb S}_{k}
=4(k+3),
\text{ and }
\dim \tilde{\mathbb S}_{2}-\dim\widetilde{\mathbb P}_{0}^3
=12.
$$
\end{proof}

Note that, given  any $\mathbf  \Phi\in \tilde{\mathbb S}_{k+2}$,
$
\overrightarrow \Delta \mathbf \Phi =- \nabla\times\nabla\times \mathbf\Phi $,
so that 
$\nabla\cdot( \overrightarrow L_k \mathbf \Phi  )= 0
$.
This also confirms the codomains of the operators of interest.

\begin{cor}\label{cor:Lap}
Given $\mathbf x_0\in\mathbb R^3$, the vector Laplacian operator restricted to the complement of homogeneous polynomial harmonic fields in the space of solenoidal fields, namely
$$
\forall k\in\mathbb N, \begin{array}{rccc}
\overrightarrow L_k|_{\tilde{\mathbb S}^*_{k+2}}:&\tilde{\mathbb S}^*_{k+2}&\to&\tilde{\mathbb S}_k,\\
&\mathbf  \Phi 
&\mapsto& \overrightarrow L_k (\mathbf  \Phi )
\end{array}
\text{ or }
\begin{array}{rrcc}
\overrightarrow L_0|_{\tilde{\mathbb S}^*_2}:&\tilde{\mathbb S}^*_2&\to&\widetilde{\mathbb P}_{0}^3,\\
&\mathbf  \Phi &\mapsto&\overrightarrow L_0\mathbf  \Phi ,
\end{array}
$$
 is also surjective.
 As a consequence, 
$$
\forall k\in\mathbb N_0, \dim\ker \left(\overrightarrow L_k|_{\tilde{\mathbb S}^*_{k+2}}\right) = 2k+5.
$$
\end{cor}
\begin{proof}
Since for all $\mathbf  \Phi\in \tilde{\mathbb S}_{k+2}$, $\overrightarrow \Delta \mathbf \Phi =- \nabla\times\nabla\times \mathbf\Phi $, then $ \widetilde{\mathbb H}_{k+2}\subset \ker \left(\overrightarrow L_k|_{\tilde{\mathbb S}_{k+2}}\right)$.
Moreover, by definition, for all $k\in\mathbb N_0$, $\widetilde{ \mathbb S}_{k+2} = \widetilde{  \mathbb S}_{k+2}^*\oplus \widetilde{\mathbb H}_{k+2}$.
Therefore, following Proposition \ref{prop:VecLapRestrict}, $\overrightarrow L_k|_{\tilde{\mathbb S}^*_{k+2}}$ is indeed surjective.

Finally, the dimension of the kernel follows from the fact that according to Proposition \ref{prop:dimHarm} for all $ k\in\mathbb N_0$, $\dim \widetilde{\mathbb H}_k = 2k+3$:
$$
\dim\ker \left(\overrightarrow L_k|_{\tilde{\mathbb S}^*_{k+2}}\right) 
=\dim\ker \left(\overrightarrow L_k|_{\tilde{\mathbb S}_{k+2}}\right)  - \dim  \widetilde{\mathbb H}_{k+2}
=2k+5.
$$
\end{proof}

\subsection{Construction of quasi-Trefftz functions}\label{subsec:construction}
Thanks to the characterization of the quasi-Trefftz spaces from Proposition \ref{prop:QTcharact} and the two surjective operators from the previous section, individual quasi-Trefftz functions can be constructed one homogeneous polynomial component at a time, as in the scalar case. 
Besides the graded structure of polynomial spaces, here the key point is that the divergence condition ensures that the right-hand side of the Laplacian condition belongs to the desired spaces.

Algorithm \ref{alg:QT construction} proposes a procedure to construct $\displaystyle \mathbf Q \in\mathbb Q\mathbb T_p$ from the order $p$, the point $\mathbf {x}_0$, the Taylor polynomial $f_\varepsilon$ of the variable coefficient $\varepsilon$, and the polynomial harmonic vector field $\mathbf {H}$.
\begin{algorithm}[ht]
    \caption{
    $\mathbf Q = {\tt QTConstr}\left(p,\mathbf {x}_0,f_\varepsilon,\mathbf F_1,\mathbf K,\mathbf {H}\right)$
    }\label{alg:QT construction}

    {\bf Input:} {$p\in\mathbb N$,  
    $\mathbf{x}_0\in\mathbb R^3$,\\
    \nonl\phantom{\bf Input:}
    $f_\varepsilon = \displaystyle\sum_{k=0}^p \varepsilon_k$ , with $\varepsilon_k \in\hompol{k}$  $\forall k\in[\![ 0,p]\!]$ and $\varepsilon_0\neq 0$, \\
    \nonl\phantom{\bf Input:}
    $\mathbf F_1\in \widetilde{\mathbb S}^*_1$, $\mathbf K= \displaystyle\sum_{k=2}^p \mathbf{K}_k    $
    with $\mathbf K_k\in\ker \overrightarrow L_k|_{\widetilde{\mathbb S}^*_{k+2}}$, $\mathbf H= \displaystyle\sum_{k=0}^p \mathbf{H}_k    $
    with $\mathbf H_k\in\homharm{k}$
    }

    {\bf Compute}  $tmp_D:=-\varepsilon_0^{-1}G_0\varepsilon_1 \cdot\mathbf H_0$. \label{pi0}
    \\
    {\bf Find} $\mathbf G_1\in\widetilde{\mathbb I}^*_{1}$ such that  $D_0\mathbf G_1 = tmp_D$. \label{tmppi0}
    
    \For{$k\ {\rm from}\ 0\ {\rm to}\ p-2$}{
        {\bf Compute} 
        $\left\{\begin{array}{l}\displaystyle tmp_D:=-\varepsilon_0^{-1} \left(  G_{k+1}\varepsilon_{k+2}\cdot\mathbf H_{0}+\sum_{k'=1}^{k+1}\Big( G_{k+1-k'}\varepsilon_{k+2-k'}\cdot(\mathbf F_{k'}+\mathbf G_{k'}+\mathbf H_{k'})
		+\varepsilon_{k+2-k'}D_{k'-1} \mathbf G_{k'} \Big) \right),\label{tmps}\\\displaystyle tmp_L:=-\sum_{k'=0}^{k}\varepsilon_{k-k'}(\mathbf F_{k'}+\mathbf G_{k'}+\mathbf H_{k'}).
		\end{array}\right.$
		
        {\bf Find} $(\mathbf G_{k+2},\mathbf F_{k+2}^{tmp})\in \widetilde{\mathbb I}^*_{k+2}\times\widetilde{\mathbb S}^*_{k+2}$ such that $\left\{\begin{array}{l}D_{k+1}\mathbf G_{k+2} = tmp_D,\\\overrightarrow L_k \mathbf F_{k+2} = tmp_L.\end{array}\right.$\label{solvek}

        {\bf Compute } $\mathbf F_{k+2}=\mathbf F_{k+2}^{tmp}+\mathbf K_{k+2}$
    }
    
    {\bf Output: } {
    $\mathbf{Q}\in (\pol{p})^3 $  where
   $\mathbf{Q} = \mathbf {H}_0 +  \displaystyle\sum_{k=1}^p \left(\mathbf{F}_k + \mathbf{G}_k + \mathbf{H}_k\right)$ 
    }
\end{algorithm}
In order to verify that this algorithm indeed constructs a quasi-Trefftz function, there are three points to justify. 
\begin{itemize}
\item First, the temporary quantities can only be computed if they depend only on known quantities.
In step \ref{tmppi0}, this is the case as $\mathbf H_0$ is known from step \ref{pi0}.
In step \ref{tmps}, this is the case since the iterations are performed for increasing values of $k$. Indeed, given $k\in[\![0,p-2]\!]$, all quantities 
$$\mathbf H_0\text{ and }  (\mathbf F_{k'},\mathbf G_{k'},\mathbf H_{k'})\in  \widetilde{ \mathbb S}^*_{k'} \times\widetilde{ \mathbb I}^*_{k'}\times  \widetilde{\mathbb H}_{k'}\text{ for all }k'\in[\![0,k+1]\!]$$ are already known from previous steps.
\item Second, finding a solution to the divergence equation is always possible 
since all operators $D_k|_{ \widetilde{\mathbb I}^*_{k+1}}$ are surjective according respectively to Lemma \ref{lmm:divk} for $k=0$ and to Proposition \ref{prop:DivRestrict} for $k>0$.
\item Third, finding a solution to the Laplace equation is always possible. Indeed, according to Corollary \ref{cor:Lap} $\overrightarrow L_k({\tilde{\mathbb S}^{{*}}_{k+2}})={\tilde{\mathbb S}_{k}}$, and the right hand side of $tmp_L$ computed at step \ref{tmps} is precisely guaranteed to belong to ${\tilde{\mathbb S}_{k}}$ from the divergence condition as $\nabla\cdot\left(tmp_L \right) = \mathbf 0$ is equivalent to
$$
 D_{0}\mathbf  Q_{1}  =-\varepsilon_0^{-1}   G_{0}\varepsilon_{1}\cdot\mathbf H_{0}
 \text{ and }
 D_{k}\mathbf  Q_{k+1}  =-\varepsilon_0^{-1} \left( \sum_{k'=0}^{k} G_{k-k'}\varepsilon_{k+1-k'}\cdot\mathbf  Q_{k'}+\sum_{k'=1}^{k}\varepsilon_{k+1-k'}D_{k'-1} \mathbf  Q_{k'}  \right).
$$
\end{itemize}

Implementing Algorithm \ref{alg:QT construction} requires the ability to solve linear equations for the $D_{k+1}|_{\tilde{\mathbb I}^*_{k+2}}$ and $\overrightarrow L_k|_{\tilde{\mathbb S}^*_{k+2}}$ operators. In practice it is possible to proceed as follows:
\begin{itemize}
    \item construct matrices for the operators $D_{k+1}|_{\tilde{\mathbb P}_{k+2}}$ and  $C_{k+1}|_{\tilde{\mathbb P}_{k+2}}$ and compute an SVD for each of them,
    \item form bases for $\tilde{\mathbb I}^*_{k+2}$ and $\tilde{\mathbb S}^*_{k+2}$ as sets of right-singular vectors associated to non-zero singular values,
    \item construct matrices for the operators $D_{k+1}|_{\tilde{\mathbb I}^*_{k+2}}$ and $\overrightarrow L_k|_{\tilde{\mathbb S}^*_{k+2}}$, 
    \item the first matrix are square and invertible according to Lemma \ref{lmm:divk} for $k=0$ and to Proposition \ref{prop:DivRestrict} for $k>0$, so linear systems can be solved thanks to any linear solver for square systems, the second matrix is rectangular, and a solution to linear systems can be found thanks to a Least-Squares approach.
\end{itemize} 
\begin{rmk}
Alternatively, to compute an individual quasi-Trefftz function without computing any SVD, step \ref{solvek} from Algorithm \ref{alg:QT construction} together with the choice of $\mathbf H_k$ can be replaced by a single step by seeking solutions to both equations in the full solenoidal and irrotational spaces.
$$
\text{Find }(\mathbf G_{k+2},\mathbf F_{k+2})\in \left(\widetilde{\mathbb P}_{k+1}\right)^3\times\widetilde{\mathbb S}_{k+2}\text{ such that }\left\{\begin{array}{l}D_{k+1}\mathbf G_{k+2} = tmp_D,\\\overrightarrow L_k \mathbf F_{k+2} = tmp_L.\end{array}\right.$$
Such solutions exist according to Lemma \ref{lmm:divk} and Proposition \ref{prop:VecLapRestrict}.
In this case $\mathbf G_{k+2}$ might have an $\widetilde{\mathbb S}_{k+2}$  component and $\mathbf F_{k+2}$ might have a harmonic component.

This can be implemented as follows.
Following the proof of Lemma \ref{lmm:divk}, given $tmp_D$ expressed in the canonical basis as $tmp_D=\sum_{|\mathbf i|=k+1}\lambda_{\mathbf i}(\mathbf X-\mathbf x_0)^{\mathbf i}$, then 
$$
\mathbf G_{k+2}:=\sum_{|\mathbf i|=k+1}\frac{\lambda_{\mathbf i}}{i_1+1}\begin{pmatrix}
(\mathbf X-\mathbf x_0)^{\mathbf i+\mathbf e_1},\mathbf 0,\mathbf 0
\end{pmatrix}
\text{ satisfies }
D_{k+1}\mathbf G_{k+2}=tmp_D.
$$
Following the proof of Proposition \ref{prop:VecLapRestrict}, this can be performed similarly for the vector Laplacian equation given $tmp_L$  expressed in
$$
\{\mathbf \Psi^{0,1,\mathbf 0}=\overrightarrow L_0 \left(\mathbf \Psi^{2,1,2\mathbf e_2}\right),
\mathbf \Psi^{0,2,\mathbf 0}=\overrightarrow L_0 \left(\mathbf \Psi^{2,2,2\mathbf e_3}\right),
\mathbf \Psi^{0,3,\mathbf 0}=\overrightarrow L_0 \left(\mathbf \Psi^{2,3,2\mathbf e_1}\right)\} \text{ basis of }{\mathbb P}_{0}^3, 
$$
$$
\{
\mathbf \Psi^{1,1,\mathbf e_3} =\overrightarrow L_1\left(  \mathbf \Psi^{3,1,3\mathbf e_3} \right)/6,\mathbf \Psi^{1,1,\mathbf e_2} =\overrightarrow L_1\left(  \mathbf \Psi^{3,1,3\mathbf e_2} \right)/6, 
$$
$$
\mathbf \Psi^{1,2,\mathbf e_1} =\overrightarrow L_1\left(  \mathbf \Psi^{3,2,3\mathbf e_1} \right)/6,\mathbf \Psi^{1,2,\mathbf e_3} =\overrightarrow L_1\left(  \mathbf \Psi^{3,2,3\mathbf e_3} \right)/6,$$
$$ \mathbf \Psi^{1,3,\mathbf e_1} =\overrightarrow L_1\left(  \mathbf \Psi^{3,3,3\mathbf e_1} \right)/6, \mathbf \Psi^{1,3,\mathbf e_2}=\overrightarrow L_1\left(  \mathbf \Psi^{3,3,3\mathbf e_2} \right)/6 ,$$
 $$\mathbf \Psi^{1,4,\mathbf e_1}=\overrightarrow L_1\left(  \mathbf \Psi^{3,4,3\mathbf e_1} \right) /6 ,  \mathbf \Psi^{1,5,\mathbf e_2}=\overrightarrow L_1\left(  \mathbf \Psi^{3,5,3\mathbf e_2} \right) /6
 \}\text{ basis of }\tilde{\mathbb S}_1, 
$$
or if $k>1$ in the basis of $\tilde{\mathbb S}_k$ defined as the union of
$$\{\overrightarrow L_k (\mathbf \Psi^{k+2,1,\mathbf i}), \forall \mathbf i\in\mathbb N_0^3\text{ such that }|\mathbf i|=k+2,i_1 = 0, i_2\leq k\},$$
$$\{\overrightarrow L_k (\mathbf \Psi^{k+2,2,\mathbf i}), \forall \mathbf i\in\mathbb N_0^3\text{ such that }|\mathbf i|=k+2,i_2 = 0, i_3\leq k\},$$
 $$\{\overrightarrow L_k (\mathbf \Psi^{k+2,3,\mathbf i}), \forall \mathbf i\in\mathbb N_0^3\text{ such that }|\mathbf i|=k+2,i_3 = 0, i_1\leq k\},$$
 $$\{\overrightarrow L_k (\mathbf \Psi^{k+2,4,\mathbf i}), \forall \mathbf i\in\mathbb N_0^3\text{ such that }|\mathbf i|=k+2,i_1>2\}, \text{ and }$$
 $$\{\overrightarrow L_k (\mathbf \Psi^{k+2,5,\mathbf i}), \forall \mathbf i\in\mathbb N_0^3\text{ such that }|\mathbf i|=k+2,i_2>2\}.$$
 \end{rmk}

{
 From Algorithm \ref{alg:QT construction}, it is then natural to construct a quasi-Trefftz basis as follows. Given any basis $\{
    \mathbf b_j 
    \}_{j\in [\![ 1, (p+1)(p+3)]\!]}$ of $\cup_{k\in[\![0,p]\!]} \homharm{k}$, any basis $\{ \mathbf s_j \}_{j\in [\![ 1, 3]\!]}$ of $\widetilde{\mathbb S}^*_1$, and any basis $\{\mathbf k_j\}_{j\in [\![ 1, (p-1)(p+3)]\!]}$ of $\cup_{k\in[\![2,p]\!]} \ker \overrightarrow L_k|_{\widetilde{\mathbb S}^*_{k+2}}$, then 
    \begin{equation*}
        \left\{ 
        {\mathbf {Q}_j = \tt QTConstr}\left(p,\bm{x}_0,T_p [\varepsilon],
        \mathbf s_j, \mathbf 0^3, \mathbf 0^3
        \right)
        \right\}_{j\in [\![ 1, 3]\!]}
        \cup
        \left\{ 
        {\mathbf {Q}_j = \tt QTConstr}\left(p,\bm{x}_0,T_p [\varepsilon], \mathbf 0^3,
        \mathbf k_j, \mathbf 0^3
        \right)
        \right\}_{j\in [\![ 1, (p-1)(p+3)]\!]}
    \end{equation*}
    $$
        \cup
        \left\{ 
        {\mathbf {Q}_j = \tt QTConstr}\left(p,\bm{x}_0,T_p [\varepsilon], \mathbf 0^3, \mathbf 0^3,
        \mathbf b_j
        \right)
        \right\}_{j\in [\![ 1, (p+1)(p+3)]\!]}
    $$
    forms a basis of $\mathbb{QT}_{p}$.
    In practice, an example of basis for $\widetilde{\mathbb S}^*_1$ is given by
    $$
        \{
        \begin{pmatrix}(\mathbf X-\mathbf x_0)^{\mathbf e_2},-(\mathbf X-\mathbf x_0)^{\mathbf e_1},\mathbf 0\end{pmatrix},
        \begin{pmatrix}(\mathbf X-\mathbf x_0)^{\mathbf e_3},\mathbf 0,-(\mathbf X-\mathbf x_0)^{\mathbf e_1}\end{pmatrix},
        \begin{pmatrix}\mathbf 0,(\mathbf X-\mathbf x_0)^{\mathbf e_3},-(\mathbf X-\mathbf x_0)^{\mathbf e_2}\end{pmatrix}
        \}
    $$
    an explicit algorithm to construct a basis for $\cup_{k\in[\![0,p]\!]} \homharm{k}$ is presented in \cite{fontana2026quasitrefftzspacesfirstordertimeharmonic}, and a basis for each $\ker \overrightarrow L_k|_{\widetilde{\mathbb S}^*_{k+2}}$ can be constructed numerically as the set of right-singular vectors associated to non-zero singular values from a computed SVD of any matrix of the operator $ \overrightarrow L_k|_{\widetilde{\mathbb S}^*_{k+2}}$.
    }

\subsection{Dimension}
One can now establish the dimension of the quasi-Trefftz spaces.
\begin{prop}
Given $p\in\mathbb N$ such that $p>2$ and $\mathbf x_0\in\mathbb R^3$,
consider Maxwell's equation \eqref{eq:Max} 
with scalar variable coefficients $\varepsilon\in\mathcal C^{p}(\mathbf x_0)$. The dimension of associated quasi-Trefftz spaces is
$$
\dim\mathbb Q\mathbb T_p=2p^2+6p+3.
$$
\end{prop}
\begin{proof}
Following Algorithm \ref{alg:QT construction}, combined with Proposition \ref{prop:dimHarm}, Corollary \ref{cor:dim*s}, Proposition \ref{prop:DivRestrict} and Proposition \ref{prop:VecLapRestrict}, the dimension of the quasi-Trefftz spaces can be computed as
$$\dim( {\mathbb P}_{0} )^3+ \dim \widetilde{\mathbb S}^*_{1}+\dim\widetilde{\mathbb H}_{1}
+\sum_{k=0}^{p-2}\left(\dim\ker \left(\overrightarrow L_k|_{\tilde{\mathbb S}^*_{k+2}}\right)+\dim\widetilde{\mathbb H}_{k+2}\right)
= 3+3+5+\sum_{k=0}^{p-2}( 2k+5 + 2k+7)
,
$$
which concludes the proof.
\end{proof}

\begin{rmk}
In \cite{CessThesis} (see Section II.9.1), a discrete PW space of dimension $(p+3)(p+1)$ is proposed to approximate solutions of a single first order equation $\nabla\times\mathbf V = \omega\mathbf V$.
This would correspond to a space of dimension $2(p+3)(p+1)$ to approximate solutions of two first order equations or equivalently a second order equation $\nabla\times\nabla\times\mathbf V = \omega^2\mathbf V$.
For comparison, naive quasi-Trefftz spaces 
$$
\left\{ 
 \mathbf  Q \in \left({\mathbb P}_{p}\right)^3,
 T_{p-2}\left[ \nabla\times\nabla\times\mathbf  Q -\varepsilon \mathbf  Q \right] = \mathbf 0^3
\right\}
=
\left\{ 
 \mathbf  Q \in \left({\mathbb P}_{p}\right)^3,
 T_{p-2}\left[ \nabla\times\nabla\times\mathbf  Q -\varepsilon \mathbf  Q \right] = \mathbf 0^3,
  T_{p-2}\left[\nabla\cdot\left(\varepsilon\mathbf  Q  \right) \right] = \mathbf 0
\right\}
$$
would precisely be of dimension 
$2(p+3)(p+1)
=2p^2+8p+6
$, while the proposed quasi-Trefftz spaces $\mathbb Q\mathbb T_p$ introduced in Definition \ref{dfn:QTp} is smaller due to the choice to impose the divergence condition to a higher order.
\end{rmk}

{
    The dimension of the quasi-Trefftz spaces $\mathbb{QT}_p$ is minimal among spaces defined by quasi-Trefftz conditions of the form
    \begin{equation*}
        \{\mathbf  Q \in \left({\mathbb P}_{p}\right)^3,
         T_{q_1}\left[ \nabla\times\nabla\times\mathbf  Q -\varepsilon \mathbf  Q \right] = \mathbf 0^3,
          T_{q_2}\left[\nabla\cdot\left(\varepsilon\mathbf  Q  \right) \right] = \mathbf 0\},
    \end{equation*}
    with $(q_1, q_2) \in (\mathbb N_0)^2$, while retaining the optimal approximation properties. Indeed,
    \begin{itemize}
        \item choosing $q_1 < p-2$ or $q_2 < p-1$, increases the the dimension of the corresponding space;
        \item choosing $q_1 > p-2$ or $q_2 > p-1$, makes the space loose the best approximation properties since $T_p[\mathbf V]$ is not one of its elements.
    \end{itemize}
    In that sense, the dimension of quasi-Trefftz spaces $\mathbb{QT}_p$ is optimal: there is no smaller quasi-Trefftz space that retains the same best approximation properties.
}

\section{Numerical results}

This section presents numerical results to validate the theory developed in this article, including the construction of a basis for $\mathbb Q\mathbb T_p$ as well as the best approximation properties of this space.

Consider two numerical test cases defined on the unit cube $\Omega = [0, 1]^3$ by a variable coefficient $\varepsilon$ for $\mathcal L_M$ and an exact solution $\mathbf{V}_{\mathrm{ex}}$ of Maxwell's equation \eqref{eq:Max}: 
\begin{enumerate}
    \item {\it (Constant coefficient test)} $\varepsilon = 100$, with exact solution $\mathbf{V}_{\rm{ex}}(\mathbf x) = e^{\mathrm{i} 10 x_3} \mathbf{e}_1$;
    \item {\it (Variable coefficient test)} $\varepsilon (\mathbf x) = 64(6-x_2)$, with exact solution 
    \begin{equation*}
        \mathbf{V}_{\mathrm{ex}}(\mathbf x) = \textrm{Ai}(4(x_2-5)) e^{\mathrm{i} 8 x_3} \mathbf{e}_1,
    \end{equation*}
    where $\textrm{Ai}$ is the so-called Airy function of the first kind, i.e., the solution of the differential equation $y'' - xy = 0$ that decays exponentially as $x\to\infty$.
\end{enumerate}
The following procedure is performed at each of 50 random points $\mathbf x_0$ in $\Omega$ for each value of $p$ from $1$ to $10$:
\begin{itemize}
    \item construct a quasi-Trefftz basis $\{\varphi_j\}_{j\in[\![1,2p^2+6p+3]\!]}$ following the procedure described in Subsection~\ref{subsec:construction},
    \item construct a quasi-Trefftz approximation $\mathbf{V}^p_{\mathrm{ap}} \in \mathbb{QT}_p$ by matching the Taylor coefficients of its expansion centered at $\mathbf{x}_0$ with those of the exact solution $\mathbf{V}_{\mathrm{ex}}$, computing the numerical solution as the Least-Squares solution,
    \item for decreasing values of $r$, estimate the error $\bm{err}(r,p, \mathbf{x}_0) \coloneqq \lVert \mathbf{V}_{\rm ex} - \mathbf{V}^p_{\text{ap}} \rVert_{\infty, \mathcal B(\mathbf{x}_0, r)}$.
\end{itemize}
For each value of $p$, consider the worst error 
$$
\bm{err}(r,p) \coloneqq \max_{\mathbf{x}_0} \bm{err}(r,p, \mathbf{x}_0)
$$
as a function of $r$; according to Theorem \ref{thn:BAP}, the order of convergence predicted by the theory is $p+1$.
Figure~\ref{fig:convergence} displays the corresponding numerical results, and for both test cases the observed orders of convergence match the predicted ones for each value of $p$.

\begin{figure}[htbp]
    \centering
    \begin{subfigure}{0.55\textwidth}
        \centering
        \resizebox{!}{0.6\textwidth}{%
        \input{test_1}
        }
        \caption{Constant coefficient test case.}
        \label{fig:convergence test_mesh_3 plots}
    \end{subfigure}
    \begin{subfigure}{0.4\textwidth}
        \centering
        \resizebox{!}{0.8\textwidth}{%
        \input{test_2}
        }
        \caption{Variable coefficient test case.}
        \label{fig:convergence test_mesh_4 plots}
    \end{subfigure}
    \\
    \begin{subfigure}{\textwidth}
        \bigskip
        \centering
        \resizebox{!}{0.14\textwidth}{%
        \begin{tabular}{|c|c|c|c|c|c|c|c|c|c|c|}
            \hline
            Radius $r$ & $p=1$ & $p=2$ & $p=3$ & $p=4$ & $p=5$ & $p=6$ & $p=7$ & $p=8$ & $p=9$ & $p=10$ \\
            \hline\hline
            $\lambda/2$ & 0.75 & 2.27 & 3.50 & 4.63 & 5.71 & 6.77 & 7.81 & 8.84 & 9.86 & 10.88 \\
            $\lambda/2^2$ & 1.69 & 2.80 & 3.86 & 4.90 & 5.92 & 6.94 & 7.95 & 8.96 & 9.96 & 10.97 \\
            $\lambda/2^3$ & 1.93 & 2.95 & 3.96 & 4.97 & 5.98 & 6.98 & 7.99 & 8.99 & 9.99 & 10.99 \\
            $\lambda/2^4$ & 1.98 & 2.99 & 3.99 & 4.99 & 5.99 & 7.00 & 8.00 & 9.00 & 10.00 & 11.00 \\
            $\lambda/2^5$ & 2.00 & 3.00 & 4.00 & 5.00 & 6.00 & 7.00 & 8.00 & 9.00 & 9.84 & - \\
            $\lambda/2^6$ & 2.00 & 3.00 & 4.00 & 5.00 & 6.00 & 7.00 & 8.00 & - & - & - \\
            $\lambda/2^7$ & 2.00 & 3.00 & 4.00 & 5.00 & 6.00 & 6.99 & - & - & - & - \\
            $\lambda/2^8$ & 2.00 & 3.00 & 4.00 & 5.00 & 6.00 & - & - & - & - & - \\
            \hline
        \end{tabular}
        }
        \caption{Constant coefficient test case.}
        \label{fig:table test_mesh_3}
    \end{subfigure}
    \\
    \begin{subfigure}{\textwidth}
        \bigskip
        \centering
        \resizebox{!}{0.14\textwidth}{%
        \begin{tabular}{|c|c|c|c|c|c|c|c|c|c|c|}
            \hline
            Radius $r$ & $p=1$ & $p=2$ & $p=3$ & $p=4$ & $p=5$ & $p=6$ & $p=7$ & $p=8$ & $p=9$ & $p=10$ \\
            \hline\hline
            $\lambda/2$ & 0.84 & 2.04 & 3.11 & 4.21 & 5.36 & 6.46 & 7.54 & 8.61 & 9.66 & 10.71 \\
            $\lambda/2^2$ & 1.43 & 2.51 & 3.62 & 4.71 & 5.78 & 6.83 & 7.86 & 8.88 & 9.90 & 10.92 \\
            $\lambda/2^3$ & 1.81 & 2.84 & 3.90 & 4.93 & 5.94 & 6.96 & 7.97 & 8.97 & 9.98 & 10.98 \\
            $\lambda/2^4$ & 1.96 & 2.96 & 3.97 & 4.98 & 5.99 & 6.99 & 7.99 & 8.99 & 10.00 & 11.00 \\
            $\lambda/2^5$ & 2.00 & 2.99 & 3.99 & 5.00 & 6.00 & 7.00 & 8.00 & 9.00 & 10.00 & 10.82 \\
            $\lambda/2^6$ & 2.01 & 3.00 & 4.00 & 5.00 & 6.00 & 7.00 & 8.00 & 8.94 & - & - \\
            $\lambda/2^7$ & 2.00 & 3.00 & 4.00 & 5.00 & 6.00 & 7.00 & 7.53 & - & - & - \\
            $\lambda/2^8$ & 2.00 & 3.00 & 4.00 & 5.00 & 6.00 & - & - & - & - & - \\
            \hline
        \end{tabular}
        }
        \caption{Variable coefficient test case.}
        \label{fig:table mesh_test_4}
    \end{subfigure}
    \caption{Error convergence for the two numerical tests with $\varepsilon = 100$ and $\varepsilon = 64(6-X_2)$, respectively: convergence plots ({\it top}) and convergence rates tables ({\it middle} and {\it bottom}).
    In the tables, $\lambda = 2\pi/\sqrt{\varepsilon} = \pi/5$ in the constant coefficient case, while $\lambda = \pi/4$ in the variable coefficient one.}
    \label{fig:convergence}
\end{figure}

\section{Conclusion}
The general framework of the work presented here is the following.
Consider a vector space $\mathbb V$, as well as an exact sequence of differential operators (each one of them being of order one)
$$\{0\}\rightarrow \mathbb V\xrightarrow{\mathcal L_1} \mathbb V^3\xrightarrow{\mathcal L_2} \mathbb V^3\xrightarrow{\mathcal L_3} \mathbb V\rightarrow \{0\}.$$
While it is more standard to study these operators defined between standard spaces of polynomials, 
this article strongly leverages the graded structure of polynomial spaces and properties of homogeneous polynomial spaces.
In the general context of polynomial spaces, the differential operator for Maxwell's equation is $\mathcal L_M =\mathcal L_2 \circ\mathcal L_2-\varepsilon Id $ where $\mathcal L_2$ is the curl operator, $\varepsilon$ is a variable coefficient and $Id$ is the identity operator. 
This operator acts between spaces of vector fields, and this 

The goal of this work is to introduce a notion of polynomial quasi-Trefftz spaces for this differential operator $\mathcal L_M$, with the same high-order best approximation properties as the corresponding standard polynomial space but with an optimal dimension, and to propose a procedure to construct a quasi-Trefftz basis. 

Back to the general framework, assuming that the operator $\mathcal L_3\circ\mathcal L_1$ is surjective leads to the existence of a Helmholtz decomposition.
Then defining $\mathbb H:=\ker \mathcal L_2\cap\ker \mathcal L_3\subset \mathbb V^3$, and assuming that $\dim\ker \mathcal L_2+\dim\ker \mathcal L_3-\dim \mathbb H=3\dim \mathbb V$
leads to the uniqueness of the decomposition.
To be more explicit, one can define ${\mathbb S}^*$ and ${\mathbb I}^*$ such that
$$
\ker \mathcal L_2 = {\mathbb I}^*\oplus\mathbb H
\text{ and }
\ker \mathcal L_3 ={\mathbb S}^*\oplus\mathbb H,
$$
and this then implies that $ \mathbb V^3 =  {\mathbb I}^*\oplus\mathbb H\oplus {\mathbb S}^*$.
Moreover, as in the scalar case, the highest-order derivative terms in $\mathcal L_M$ play a crucial role in this work.
Here the component-wise version of the operator $\mathcal L_3\circ\mathcal L_1$, denoted $\overrightarrow{\mathcal L_3\circ\mathcal L_1}$, comes into play to leverage the unique decomposition assuming that $\mathcal L_2 \circ\mathcal L_2=-\overrightarrow{\mathcal L_3\circ\mathcal L_1}+\mathcal L_1\circ\mathcal L_3$.
However, unlike in the scalar case, differential operators between spaces of polynomial vector fields are not necessarily surjective, and in particular $\mathcal L_2 \circ\mathcal L_2$ is not surjective.
The well-documented divergence condition is a consequence of the exact sequence definition: $\mathcal R(\mathcal L_2)=\ker \mathcal L_3$, therefore $\mathcal L_M\varphi=0$ implies that $\varepsilon\varphi\in\ker\mathcal L_3$.
This fact, together with the unique decomposition, can then be used to express $\mathcal L_M\varphi=0$ equivalently as the combination of a divergence-free condition for the full vector field $\varphi$ and a differential equation in which the only second order term is the $\overrightarrow{\mathcal L_3\circ\mathcal L_1}$ operator acting only on the $\mathbb S^*$ component.
This combined expression now shares with the scalar case the fact that the highest-order derivative operators, acting between well-chosen spaces, are surjective.

This fact, leveraged in the quasi-Trefftz setting, is precisely the cornerstone of the present study of the quasi-Trefftz spaces $\mathbb Q\mathbb T_p$ for Maxwell's equation $\mathcal L_M\varphi=0$.

\bibliographystyle{alpha}
\bibliography{biblio}

%
%
%
%
%

\end{document}